\documentclass[a4paper]{article}
\usepackage[top=25truemm,bottom=25truemm,left=20truemm,right=20truemm]{geometry}
\usepackage{amsmath,amsthm,amsfonts,amssymb,stmaryrd,indentfirst,tabularray,makecell,hhline,comment,mathtools}
\usepackage[mathscr]{euscript}
\usepackage[shortlabels]{enumitem}
\usepackage[labelfont=rm]{caption}
\usepackage{titlesec}
\usepackage{tikz-cd}
\usepackage{url}
\usepackage[backref]{hyperref}
\hypersetup{
    colorlinks = true,
    linkcolor = [rgb]{0.3,0.55,0.25},
    citecolor = [rgb]{0.164,0.39,0.77},
    urlcolor = [rgb]{0.7,0.2,0.4}
}

\DeclareMathOperator{\id}{id}
\DeclareMathOperator{\Ker}{Ker}
\DeclareMathOperator{\Imag}{Im}
\DeclareMathOperator{\Der}{Der}

\DeclareMathOperator{\tDer}{tDer}

\DeclareMathOperator{\End}{End}
\DeclareMathOperator{\Aut}{Aut}
\DeclareMathOperator{\tAut}{tAut}
\DeclareMathOperator{\Tr}{Tr}
\DeclareMathOperator{\Div}{Div}

\DeclareMathOperator{\sdiv}{div} 
 
\DeclareMathOperator{\DR}{DR}

\DeclareMathOperator{\Ad}{Ad}
\DeclareMathOperator{\ad}{ad}
\DeclareMathOperator{\Exp}{exp}

\DeclareMathOperator{\rot}{rot}
\DeclareMathOperator{\irr}{irr}

\newcommand{\K}[1]{\mathbb{K}[[#1]]}

\newcommand{\op}[0]{^\mathrm{op}}
\newcommand{\argdot}[0]{\,\cdot\,}
\newcommand{\omid}[0]{\langle\omega\rangle}
\newcommand{\gr}[0]{_\mathrm{gr}}

\theoremstyle{plain}
\newtheorem{theorem}{Theorem}
\newtheorem*{theorem*}{Theorem}
\newtheorem{proposition}[theorem]{Proposition}

\newtheorem{lemma}[theorem]{Lemma}
\newtheorem{corollary}[theorem]{Corollary}
\newtheorem{conjecture}[theorem]{Conjecture}

\theoremstyle{definition}
\newtheorem{definition}[theorem]{Definition}
\newtheorem{remark}[theorem]{Remark}
\newtheorem{example}[theorem]{Example}
\newtheorem{computation}[theorem]{Computation}

\numberwithin{theorem}{section}
\everymath{\displaystyle}
\allowdisplaybreaks

\titleformat{\section}{\large\scshape}{\thesection.}{3pt}{}
\titleformat{\subsection}{\scshape}{\thesubsection}{3pt}{}

\usepackage{titling}  
\pretitle{
	\begin{center}\Large\bfseries
}
\posttitle{
	\end{center}\ \\[-20pt]
}
\preauthor{
	\begin{center}\large
}
\postauthor{
	\end{center}\ \\[-60pt]
}

\title{The Formality of the Goldman--Turaev Lie Bialgebra\\ on a Closed Surface}
\author{Toyo TANIGUCHI \thanks{Graduate School of Mathematical Sciences, The University of Tokyo. 3-8-1, Komaba, Meguro-ku, Tokyo, 153-8914, Japan. E-mail: \texttt{toyo(at)ms.u-tokyo.ac.jp}}}
\date{}

\begin{document}
\maketitle

\begin{abstract}
\noindent We reformulate the Kashiwara--Vergne groups and associators in higher genera, introduced in Alekseev--Kawazumi--Kuno--Naef, in terms of non-commutative connections using the tools developed in a previous paper. As the main result, the case of closed surfaces is dealt with to determine the pro-unipotent automorphism group of the associated graded of the Goldman--Turaev Lie bialgebra.\\
\end{abstract}

\noindent{\textit{2020 Mathematics Subject Classification: 16D20, 17A61, 53D30, 57K20, 58B34. }

\section{Introduction}
The Kashiwara--Vergne (KV) problem originated from Lie theory asking the existence of certain infinite Lie series, the \textit{Kashiwara--Vergne associators}. The problem is partially solved in their original paper \cite{kv} and then by Rouvi\`ere \cite{rouviere} and Vergne \cite{vergne}, and completely solved by Alekseev--Meinrenken \cite{am}. Later, by a remarkable result by Alekseev--Torossian \cite{at}, the KV problem is reduced to the existence of Drinfeld associators, which is known over the field of characteristic zero in the original paper \cite{drinfeld} by Drinfeld. The space of solutions to the KV problem is a bi-torsor over the groups KV and KRV, which are called the Kashiwara--Vergne groups.

The defining relations of the KV associators and the KV groups have a higher genus analogue. Let $\mathbb{K}$ be a field of characteristic zero and $\pi = \pi_1(\Sigma)$ the fundamental group of a compact oriented connected surface $\Sigma$. Then, the \textit{trace space} $|\mathbb{K}\pi| = \mathbb{K}\pi/[\mathbb{K}\pi,\mathbb{K}\pi]$ is endowed with the structure of the Goldman--Turaev Lie bialgebra. Alekseev--Kawazumi--Kuno--Naef \cite{akkn} shows that, in the genus $0$ case, the formality problem of the Goldman--Turaev Lie bialgebra, which asks if the completion of this Lie bialgebra under a suitable filtration is isomorphic to its associated graded (see Section \ref{sec:kvboundary} for the detail), is essentially equivalent to the KV problem; this perspective immediately allows us to define the higher genus KV associators as the solutions to the formality problem and the KV groups as the (pro-unipotent) automorphism group of the (graded) Goldman--Turaev Lie bialgebra. We remark that, also in the genus $0$ case, a direct construction of a formality isomorphism from a Drinfeld associator by means of the LMO functor is originally given in \cite{mas} by Massuyeau. On the other hand, we have higher genus analogues of the Drinfeld associators. For genus $1$, they are called \textit{elliptic associators} introduced by Enriquez \cite{enriquez}, and their relation to the corresponding KV problem is discussed in \cite{akkn}. For an arbitrary genus, several generalisations of the GT groups and the Drinfeld associators are proposed: \cite{gonzalez} by Gonzalez, \cite{felder} by Felder, and Section 5 of \cite{ciw} by Campos--Idrissi--Willwacher. They all agree in genus $1$, but the relation between them and the KV associators for genus $\geq 2$ is still an open question.

Getting back to the KV problem in higher genera, the existence problem itself is completely solved in \cite{akkn}, and the key step in the proof is the factorisation of the Turaev cobracket in two parts: a Hamiltonian flow and a non-commutative divergence. In this paper, we solve the (non-)uniqueness part of the formality problem on a closed surface: we determine the pro-unipotent part of the automorphism group of the Goldman--Turaev Lie bialgebra on a closed surface using the factorisation of the Turaev cobracket obtained in the author's previous paper \cite{toyo}. We start with a reformulation of the KV groups and the space of solution $\mathrm{SolKV}$ for a surface $\Sigma$ with non-empty boundary in terms of connections in non-commutative geometry and later see the case of closed surfaces based on the observation.

Let $H = H_1(\Sigma;\mathbb{K})$ be the first homology group of the surface $\Sigma$. In their formulation, the KV groups are described in terms of group $1$-cocycles on the group $\tAut_\mathrm{Hopf}(\widehat{\mathbb{K}\pi})$ of \textit{tangential automorphisms} on the complete Hopf algebra $\widehat{\mathbb{K}\pi}$, and the graded counterpart $\tAut_\mathrm{Hopf}(\hat T(H))$ on the completed tensor algebra $\hat T(H)$. We first see that a tangential automorphism of $\widehat{\mathbb{K}\pi}$ is interpreted as a derivation of the complete \textit{Hopf groupoid} (in the sense of Fresse \cite{fresse}) $\widehat{\mathbb{K}\mathscr{G}}$ associated with the fundamental groupoid $\mathscr{G}$ of $\Sigma$, and similarly for the graded counterpart $\hat{\mathscr{T}}$. We also have their subcategories $\partial\widehat{\mathbb{K}\mathscr{G}}$ and $\partial\hat{\mathscr{T}}$ generated by \textit{boundary loops}.

Next, we define connections $ \nabla'\!_{\mathcal{C},\mathsf{fr}}$ and $\nabla'\!\!_{H,\mathsf{fr}}$ on some $\mathbb{K}\mathscr{G}$-module and $\hat{\mathscr{T}}$-module, respectively, depending on a framing $\mathsf{fr}$ on $\Sigma$ and a free-generating system $\mathcal{C}$ of $\pi$. Recall that a connection on an $A$-module $N$, where $A$ is a $\mathbb{K}$-algebra, is a $\mathbb{K}$-linear map $\nabla\colon N\to \Omega^1A\otimes_AN$ satisfying the Leibniz rule, where $\Omega^1A$ is the space of non-commutative $1$-forms. Then, the value of those $1$-cocycles on an automorphism $G$ is written as (the trace of) the difference of connections by the action of $G$: $\Tr(\Ad_G\!\nabla - \nabla)$ where $\nabla$ is either of the two above. This takes its value in the space of \textit{de Rham 1-forms} $\DR^1\!A$, a quotient space of the space $\Omega^1A$ of non-commutative $1$-forms. We have the exterior derivative $d$, which sends an element of the trace space $|A| = A/[A, A]$ to $\DR^1\!A$. One of the main results is the following.

\begin{theorem*}[Theorem \ref{thm:kvboundary}]
The KV groups and the set of KV associators are expressed as follows:
\begin{itemize}
	\item $\mathrm{KV}^\mathsf{fr}_{(g,n+1)} = \Big\{G \in\Aut^+_{\mathrm{Hopf},\partial}(\widehat{\mathbb{K}\mathscr{G}}): \Tr(\Ad_G\!\nabla'\!_{\mathcal{C},\mathsf{fr}} - \nabla'\!_{\mathcal{C},\mathsf{fr}}) \in d|\partial\widehat{\mathbb{K}\mathscr{G}}|\Big\}$,
	\item $\mathrm{KRV}^\mathsf{fr}_{(g,n+1)} = \Big\{G \in\Aut^+_{\mathrm{Hopf},\partial}(\hat{\mathscr{T}}):\Tr(\Ad_G\!\nabla'\!\!_{H,\mathsf{fr}} - \nabla'\!\!_{H,\mathsf{fr}}) \in d|\partial\hat{\mathscr{T}}|\Big\}$, and
	\item $\mathrm{SolKV}^\mathsf{fr}_{(g,n+1)} = \Big\{G\in\mathrm{Isom}^+_{\mathrm{Hopf},\partial}(\hat{\mathscr{T}},\widehat{\mathbb{K}\mathscr{G}}): \Tr(\Ad_G\!\nabla'\!\!_{H,\mathsf{fr}} - \nabla'\!_{\mathcal{C},\mathsf{fr}}) \in d|\partial\widehat{\mathbb{K}\mathscr{G}}|\Big\}$.
\end{itemize}
\end{theorem*}
\noindent Here, the superscript $^+$ indicates the pro-unipotent part of each space. In particular, the Lie algebra corresponding to $\mathrm{KRV}^\mathsf{fr}_{(g,n+1)}$ is given by
\[
	\mathfrak{krv}^\mathsf{fr}_{(g,n+1)} = \Big\{g \in\Der^+_{\mathrm{Hopf},\partial}(\hat{\mathscr{T}}):\sdiv^{\nabla'\!\!_{H,\mathsf{fr}}}(g) \in |\partial\hat{\mathscr{T}}|\Big\}.
\]
The divergence map $\sdiv$ is defined from a connection, whose construction is given in \cite{toyo}.

Now consider the case of closed surfaces. The existence of a formality morphism follows from the case of a surface with only one boundary component, which is solved in \cite{akkn}, so it is enough to consider the automorphism group of the associated graded. Based on the reformulation above, we have a similar description in the case of closed surfaces:

\begin{theorem*}[Theorem \ref{thm:kvclosed}]
The pro-unipotent part of the automorphism group of the associated graded of the Goldman--Turaev Lie bialgebra is given by $\mathrm{KRV}_{(g,0)} = \Exp(\mathfrak{krv}_{(g,0)})$, where
\[
	\mathfrak{krv}_{(g,0)} := \{g\in\Der^+(\hat L(H)_\omega)\colon \sdiv^{\nabla'_{\bullet,H}}(g) \in \Ker(|\bar\Delta_\omega|)\}\,.
\]
Here, $\hat L(H)_\omega$ is the quotient of the completed free Lie algebra over $H$ by the element $\omega\in H^{\otimes 2}$ representing the symplectic structure of $H$, $\nabla'_{\bullet,H}$ is a connection on a $\hat L(H)_\omega$-module and $\bar\Delta_\omega\colon \hat T(H)_\omega \to \hat T(H)_\omega\hat \otimes \hat T(H)_\omega$ is the reduced coproduct: $\bar\Delta_\omega(x) = \Delta_\omega(x) - x\otimes 1 - 1\otimes x$.
\end{theorem*}

The space $\Ker(|\bar\Delta_\omega|)$ deserves some description, so we compute it in small degrees at the end of the paper.\\

\noindent\textbf{Organisation of the paper.} In Sections \ref{sec:gtlb} and \ref{sec:kvboundary}, we recall the definition of the Goldman--Turaev Lie bialgebra, the formality problem, and main results in \cite{akkn}. Section \ref{sec:reform} is occupied with the reformulation of the KV groups. In Section \ref{sec:rewriting}, we construct a basis of $|T(H)_\omega|$ for later sections. The formality problem on a closed surface is dealt with in Section \ref{sec:closed}. Section \ref{sec:redcoprod} is devoted to the computation of the kernel of the reduced coproducts.\\

\noindent\textbf{Acknowledgements.} The author thanks Yusuke Kuno for a discussion on Section \ref{sec:redcoprod}, Nariya Kawazumi for thoroughly reading the draft and Geoffrey Powell for many comments and improvements on the draft.\\

\noindent\textbf{Conventions.} $\mathbb{K}$ is a field of characteristic zero. Unadorned tensor products are always over $\mathbb{K}$.\\

\section{The Goldman--Turaev Lie Bialgebra}\label{sec:gtlb}
In this section, we recall the Goldman--Turaev Lie bialgebra and its associated graded following \cite{akkn}. Let $g\geq 0$, $n\geq -1$ and $\Sigma = \Sigma_{g,n+1}$ a connected compact oriented surface of genus $g$ and $n+1$ boundary components $\partial_0\Sigma,\dotsc, \partial_n\Sigma$. Consider the fundamental group $\pi = \pi_1(\Sigma,*)$ of $\Sigma$ and the group algebra $\mathbb{K}\pi$. We require the base point $*$ to be taken on $\partial_0\Sigma$ if $n\geq 0$. For a $\mathbb{K}$-algebra $A$, we denote the trace space $A/[A,A]$ by $|A|$. Then, $|\mathbb{K}\pi| = \mathbb{K}\pi/[\mathbb{K}\pi,\mathbb{K}\pi]$ is identified with the $\mathbb{K}$-vector space spanned by the homotopy set of free loops on $\Sigma$.

\begin{definition}
We define the \textit{Goldman bracket} $[\cdot,\cdot]\colon |\mathbb{K}\pi|^{\otimes 2} \to |\mathbb{K}\pi|$ by, for generically immersed free loops $\alpha$ and $\beta$ on $\Sigma$,
\[
	[\alpha,\beta] = \sum_{p\in\alpha\cap\beta} \mathrm{sign}(p;\alpha,\beta) \alpha*_p\beta,
\]
where $\mathrm{sign}(p;\alpha,\beta)$ is the local intersection number of the velocity vectors $\dot\alpha_p$ and $\dot\beta_p$ at $p$ with respect to the orientation of the surface, and $\alpha*_p\beta$ is the free loop obtained by first traversing $\alpha$ starting at $p$ then $\beta$.
\end{definition}

\begin{definition}
Assume $n\geq 0$. Let $\alpha$ be a free loop and  $\beta$ a based loop on $\Sigma$, both generically immersed. We define the \textit{Kawazumi--Kuno action} $\sigma\colon |\mathbb{K}\pi| \to \End_\mathbb{K}(\mathbb{K}\pi)$ by
\[
	\sigma(\alpha)(\beta) = \sum_{p\in\alpha\cap\beta} \mathrm{sign}(p;\alpha,\beta) \alpha*_p\beta.
\]
\end{definition}

\begin{definition}
We define the \textit{Turaev cobracket} $\delta\colon |\mathbb{K}\pi/\mathbb{K}1| \to |\mathbb{K}\pi/\mathbb{K}1|^{\otimes 2}$ by, for a generically immersed free loop $\alpha\colon [0,1]/\{0,1\} \to \Sigma$,
\begin{align*}
	\delta(\alpha) = \sum_{\substack{t_1\neq t_2\in[0,1]\\\alpha(t_1)=\alpha(t_2)}}  \mathrm{sign}(\alpha;t_1,t_2)\,\alpha|_{[t_1,t_2]}\otimes\alpha|_{[t_2,t_1]}\,,
\end{align*}
where $\mathrm{sign}(\alpha;t_1,t_2)$ is the local intersection number of the velocity vectors $\dot\alpha(t_1)$ and $\dot\alpha(t_2)$ with respect to the orientation of $\Sigma$.

If $\Sigma$ admits a framing $\mathsf{fr}$ (i.e., a smooth non-vanishing vector field), we define the framed version $\delta^\mathsf{fr}\colon|\mathbb{K}\pi| \to |\mathbb{K}\pi|^{\otimes 2}$ by the same formula but taking a rotation-free representative of $\alpha$. 
\end{definition}

\begin{theorem}[\cite{goldman},\cite{turaev},\cite{akkn}]
The $\mathbb{K}$-linear maps $[\cdot,\cdot]$, $\delta$ and $\delta^\mathsf{fr}$ are well-defined. The triples $(|\mathbb{K}\pi|,[\cdot,\cdot],\delta^\mathsf{fr})$ and $(|\mathbb{K}\pi/\mathbb{K}1|,[\cdot,\cdot],\delta)$ are Lie bialgebras.
\end{theorem}

\begin{theorem}[\cite{sigma}]
The map $\sigma$ is well-defined and takes its values in the space of derivations: $\sigma \colon |\mathbb{K}\pi| \to \Der_\mathbb{K}(\mathbb{K}\pi)$.
\end{theorem}

Next, we recall a weight filtration on $\mathbb{K}\pi$. Note that a (decreasing) filtration on a vector space can be pulled back along any linear map and pushed out along a surjection.
\begin{definition} Assume $n\geq 0$.
\begin{itemize}
	\item Let $\mathcal{C} = (\alpha_i,\beta_i,\gamma_j)_{1\leq i\leq g,1\leq j\leq n}$ be a free-generating system of $\pi$ so that $\alpha_i$ and $\beta_i$ form a genus pair, $\gamma_j$ is a boundary loop representing $\partial_j\Sigma$ and
	\[
		\alpha_1\beta_1\alpha_1^{-1}\beta_1^{-1}\cdots \alpha_g\beta_g\alpha_g^{-1}\beta_g^{-1}\gamma_1\cdots\gamma_n
	\]
	represents a boundary loop based at $*\in\partial_0\Sigma$ (see Figure 2 of \cite{toyo2}). Denote by $(x_i,y_i,z_j)_{1\leq i\leq g,1\leq j\leq n}$ the corresponding basis of $H = H_1(\Sigma;\mathbb{K})$.
	\item $\hat T(H) \cong \mathbb{K}\langle\!\langle x_i,y_i,z_j \rangle\!\rangle_{1\leq i\leq g, 1\leq j\leq n}$ is the completed free associative algebra over $H$ with respect to the \textit{weight grading} on $H$ defined by $\mathrm{wt}(x_i) = \mathrm{wt}(y_i) = 1$ and $\mathrm{wt}(z_j) = 2$.
	\item Consider the morphism of $\mathbb{K}$-algebras
	\[
		\theta_\mathrm{exp}\colon\mathbb{K}\pi \to \hat T(H)\colon \alpha_i \mapsto e^{x_i}, \beta_i \mapsto e^{y_i} \textrm{ and } \gamma_j \mapsto e^{z_j}.
	\]
	Define a filtration on $\mathbb{K}\pi$ by the pull-back of the weight filtration by $\theta_\mathrm{exp}$. This induces an isomorphism of Hopf algebras on the completion $\widehat{\mathbb{K}\pi}$.
	\item Define filtrations on $|\mathbb{K}\pi|$ and $|\mathbb{K}\pi/\mathbb{K}1|$ as the push-out by the natural surjections $\mathbb{K}\pi \twoheadrightarrow |\mathbb{K}\pi|\twoheadrightarrow |\mathbb{K}\pi/\mathbb{K}1|$.
	\item The filtration on $\Der_\mathbb{K}(\mathbb{K}\pi)$ is the induced one: denoting the filtration degree $\geq k$ part of $\mathbb{K}\pi$ by $F^k\mathbb{K}\pi$, $u \colon \mathbb{K}\pi \to \mathbb{K}\pi$ has filtration degree $\geq p$ if and only if  $u(F^k\mathbb{K}\pi) \subset F^{k+p}\mathbb{K}\pi$ for all $k$.
	\end{itemize}
\end{definition}

\begin{definition}
For $n=-1$, the filtration on $|\mathbb{K}\pi_1(\Sigma_{g,0})|$ is defined by the push-out along $|\mathbb{K}\pi_1(\Sigma_{g,1}) |\twoheadrightarrow |\mathbb{K}\pi_1(\Sigma_{g,0})|$ by capping the unique boundary component by a disk. The rest is analogously defined.
\end{definition}

\begin{proposition}\label{prop:filtdeg}
The $\mathbb{K}$-linear maps $[\cdot,\cdot]$, $\sigma$, $\delta$ and $\delta^\mathsf{fr}$ are compatible with the filtrations. They all have the filtration degree $(-2)$.
\end{proposition}
\noindent For the proof, see Section 3 of \cite{akkn} and references there. Following them, we denote their associated gradeds by $[\cdot,\cdot]\gr\colon |\hat T(H)|^{\otimes 2} \to |\hat T(H)|$, $\sigma\gr\colon |\hat T(H)| \to \Der_\mathbb{K}(\hat T(H))$, and so on.\\

\section{The Formality and the Kashiwara--Vergne Problem}\label{sec:kvboundary}

The \textit{formality problem} for the Goldman--Turaev Lie algebra asks if it is isomorphic to its associated graded as a Lie bialgebra and, if so, to determine the set of such isomorphisms. More precisely, a solution to the formality problem is a completed Hopf algebra isomorphism
\[
	\theta \colon \widehat{\mathbb{K}\pi} \to \mathrm{gr}(\widehat{\mathbb{K}\pi})
\]
such that $\mathrm{gr}(\theta) = \id$ and induces an isomorphism of Lie bialgebras $|\widehat{\mathbb{K}\pi}| \to |\mathrm{gr}(\widehat{\mathbb{K}\pi})|$.

In the case of a surface with non-empty boundary, it is completely solved in \cite{akkn}: the set SolKV of such isomorphisms is characterised by two equations (KVI) and (KVII). In the special case of $\Sigma = \Sigma_{0,3}$, the problem is, surprisingly, equivalent to the Kashiwara--Vergne problem in Lie theory, and the automorphism group of the Lie bialgebra and that of the associated graded are given by the KV and KRV groups. This section summarises definitions of these spaces and results in \cite{akkn}.\\

First, let us recall the definition of some spaces we need. Fix a connected oriented compact surface $\Sigma$ of genus $g$ with $(n+1)$ boundary components with $n\geq 0$ equipped with a framing $\mathsf{fr}$. Take the free-generating system $(\alpha_i,\beta_i,\gamma_j)$ of $\pi$ and the basis $(x_i,y_i,z_j)$ of $H$ as in the last section. We have the following spaces:

\begin{itemize}
	\item $\hat L(H) = L(\!(x_i,y_i,z_j)\!)_{1\leq i\leq g, 1\leq j\leq n}$ is the completed free Lie algebra over $H$, so that $\hat T(H)$ is identified with the (completed) universal enveloping algebra $U\!\hat L(H)$,
	\item $\Der^+(\hat L(H)) = \{u\colon \hat L(H)\to \hat L(H): \textrm{a continuous Lie algebra derivation, degree}\geq 1\}$,
	\item $\Aut^+(\hat L(H)) = \Exp(\Der^+(\hat L(H)))$,
	\item the space of \textit{tangential derivations}: $$\tDer^+(\hat L(H)) = \{\tilde u=(u;u_1,\dotsc,u_n): u\in \Der^+(\hat L(H)), \,u_i\in \hat L(H), \,u(z_j) = [z_j,u_j]\},$$
	\item  the space of \textit{tangential automorphisms}:
	\begin{align*}
		\tAut^+(\hat L(H)) &= \Exp(\tDer^+(\hat L(H)))\\
		& = \{\tilde G=(G;g_1,\dotsc,g_n):G\in \Aut^+(\hat L(H)),\, g_j\in\Exp(\hat L(H)), \,G(z_j) = g_j^{-1}z_jg_j\},
	\end{align*}
\end{itemize}
and some elements:
\begin{itemize}
	\item $\xi = \log\left(\prod_i(e^{x_i}e^{y_i}e^{-x_i}e^{-y_i})\prod_j e^{z_j}\right)\in \hat L(H)$,
	\item $\omega = \sum_i [x_i,y_i] + \sum_j z_j\in \hat L(H)$,
	\item $r(s) = \log\left(\frac{e^s-1}{s}\right)\in s\K s$,
	\item $\mathbf{r} = \sum_i |r(x_i) + r(y_i)|,\;\mathbf{r}' = \sum_i |r(x_i) + r(y_i)|+ \sum_j |r(z_j)|\in |\hat T(H)|$,
	\item $a_i = \rot^\mathsf{fr}(\alpha_i),\,b_i = \rot^\mathsf{fr}(\beta_i),\,c_j = \rot^\mathsf{fr}(\gamma_j) \in \mathbb{Z}$, and $\mathbf{p}^\mathsf{fr} = \sum_i |a_iy_i - b_ix_i|\in |\hat T(H)|$.
\end{itemize}

Now recall the standard divergence for a free Lie algebra and the integration of a $1$-cocycle.
\begin{definition}\ 
\begin{itemize}
	\item $\hat T(H)$ is regarded as an $\hat L(H)$-module by the left multiplication. For $w = x_i,y_i,z_j$, $d_w\colon \hat L(H) \to \hat T(H)$ is a continuous Lie algebra $1$-cocycle specified by $d_w(w') = \delta_{ww'}$ for $w' = x_i,y_i,z_j$ using Kronecker's delta.
	\item We define the \textit{single divergence} $\sdiv_{x,y,z} \colon \Der(\hat L(H)) \to |\hat T(H)|$ by 
	\[
		\sdiv_{x,y,z}(u) = \sum_{w = x_i,y_i,z_j} |d_wu(w)|.
	\]
	This is extended to $\tDer(\hat L(H))$ by the composition
	\[
		\sdiv_{x,y,z}\colon \tDer^+(\hat L(H)) \to \Der(\hat L(H)) \xrightarrow{\sdiv_{x,y,z}} |\hat T(H)|.
	\]
	The single divergence is itself a Lie algebra $1$-cocycle.
	\item For a pro-nilpotent Lie algebra $\mathfrak g$, a continuous $\mathfrak g$-module $V$ and a $1$-cocycle $\psi\colon \mathfrak g \to V$, its \textit{integration} is a group $1$-cocycle $\Psi\colon \exp(\mathfrak g) \to V$ given by, for $u\in \mathfrak g$,
	\[
		\Psi(e^u) = \frac{e^u-1}{u}\cdot \psi(u).
	\]
	The correspondence $\psi \mapsto \Psi$ is $\mathbb{K}$-linear. For the details, see Appendix A of \cite{akkn}.
	\item Since $\tDer^+(\hat L(H))$ is pro-nilpotent, we denote the integration of $\sdiv_{x,y,z}\colon \tDer^+(\hat L(H)) \to |\hat T(H)|$ by $$\mathsf{j}_{x,y,z}\colon \tAut^+(L) \to |\hat T(H)|.$$
\end{itemize}
\end{definition}

\noindent We also have many $1$-cocycles:
\begin{itemize}
	\item $\mathsf{b}^\mathsf{fr}\colon\tDer^+(\hat L(H))\to |\hat T(H)|:\tilde u\;\mapsto\;\sum_j c_j|u_j|$ and $\mathsf{c}^\mathsf{fr}\colon\tAut^+(\hat L(H))\to |\hat T(H)|$ its integration, \\[-10pt]
	\item $\sdiv^\mathsf{fr}\colon \tDer^+(\hat L(H))\to |\hat T(H)|: \tilde u\;\mapsto\; \sdiv_{x,y,z}(u)- \mathsf{b}^\mathsf{fr}(\tilde u) + u(\mathbf{r} - \mathbf{p}^\mathsf{fr})$ and $\mathsf{j}^\mathsf{fr}\colon \tAut^+(\hat L(H))\to |\hat T(H)|$ its integration, and
	\item $\sdiv^\mathsf{fr}_\mathrm{gr}\colon \tDer^+(\hat L(H))\to |\hat T(H)|: \tilde u\;\mapsto\; \sdiv_{x,y,z}(u)- \mathsf{b}^\mathsf{fr}(\tilde u)$ and $\mathsf{j}^\mathsf{fr}_\mathrm{gr}\colon \tAut^+(\hat L(H))\to |\hat T(H)|$ its integration.\\
\end{itemize}

Finally, we recall the definition of the KV groups and associators.
\begin{definition} For $g,n\geq 0$ and a framing $\mathsf{fr}$, the Kashiwara--Vergne group $\mathrm{KV}^\mathsf{fr}_{(g,n+1)}$, the graded version $\mathrm{KRV}^\mathsf{fr}_{(g,n+1)}$, and the set of the Kashiwara--Vergne associators $\mathrm{SolKV}^\mathsf{fr}_{(g,n+1)}$ is defined by the followings:
\begin{itemize}
	\item $\mathrm{KV}^\mathsf{fr}_{(g,n+1)} = \Big\{\tilde G \in\tAut^+(\hat L(H)): G(\xi) = \xi,\,\mathsf{j}^\mathsf{fr}(\tilde G) \in \Big|\sum_j z_j\K{z_j} + \xi^2\K\xi\Big|\Big\}$,
	\item $\mathrm{KRV}^\mathsf{fr}_{(g,n+1)} = \Big\{\tilde G \in\tAut^+(\hat L(H)): G(\omega) = \omega,\,\mathsf{j}^\mathsf{fr}_\mathrm{gr}(\tilde G) \in \Big|\sum_j z_j\K{z_j} + \omega^2\K\omega\Big|\Big\}$, and
	\item $\mathrm{SolKV}^\mathsf{fr}_{(g,n+1)} = \{\tilde G\in\tAut^+(\hat L(H)): G(\omega)=\xi,\,\mathsf{j}^\mathsf{fr}_\mathrm{gr}(\tilde G) -\mathbf{r} + \mathbf{p}^\mathsf{fr} \in \Big|\sum_j z_j\K{z_j} + \xi^2\K\xi\Big|\Big\}$.
\end{itemize}
\end{definition}

The groups $\mathrm{KV}^\mathsf{fr}_{(g,n+1)}$ and $\mathrm{KRV}^\mathsf{fr}_{(g,n+1)}$ are pro-unipotent since $\tAut^+(\hat L(H))$ is. The corresponding Lie algebras are given by the following:
\begin{itemize}
	\item $\mathfrak{kv}^\mathsf{fr}_{(g,n+1)} = \Big\{\tilde g \in\tDer^+(\hat L(H)): g(\xi) = 0,\,\sdiv^\mathsf{fr}(\tilde g) \in \Big|\sum_j z_j\K{z_j} + \xi^2\K\xi\Big|\Big\}$, and
	\item $\mathfrak{krv}^\mathsf{fr}_{(g,n+1)} = \Big\{\tilde g \in\tDer^+(\hat L(H)): g(\omega) = 0,\,\sdiv^\mathsf{fr}_\mathrm{gr}(\tilde g) \in \Big|\sum_j z_j\K{z_j} + \omega^2\K\omega\Big|\Big\}$.
\end{itemize}

One of the main results in \cite{akkn} is the following.

\begin{theorem}[\cite{akkn}, Theorem 6.27] For $\Sigma = \Sigma_{g,n+1}$ with $n\geq 0$, an isomorphism of filtered Hopf algebras $\theta\colon\widehat{\mathbb{K}\pi} \to \hat T(H)$ with $\mathrm{gr}(\theta) = \id$ gives a solution to the formality problem if and only if $\theta\circ\theta_\mathrm{exp}^{-1}$ lifts to an element in $\mathrm{SolKV}_{(g,n+1)}$ up to conjugation by an element of $\Aut^+(\hat L(H))$.
\end{theorem}

SolKV is a (two-sided) torsor over KV and KRV, which is apparent from their defining equations.\\

\section{Reformulation of KV Groups}\label{sec:reform}
In this section, we reformulate the defining equations of the KV and KRV groups using certain connections in non-commutative geometry. We use notations in \cite{toyo} and \cite{toyo2}. In Section \ref{sec:closed}, we will deal with the case of closed surfaces based on observations made here.

\subsection{Tangential Derivations}
We start with tangential derivations, which are best understood in the setting of linear categories. Let $g,n\geq 0$ and denote $(n+1)$ boundary components of $\Sigma_{g,n+1}$ by $\partial_0\Sigma,\dotsc,\partial_n\Sigma$. Take one base point $*_j$ for each $j$ and put  $V = \{*_j\}_{0\leq j\leq n}$. Now set $\mathscr{G} = \pi_1(\Sigma_{g,n+1},V)$, the fundamental groupoid of $\Sigma_{g,n+1}$ with base points $V$.

\begin{definition}
Let $(\gamma_j)_j$ be as above: $\gamma_j$ is a simple loop based at $*_0$ representing the $j$-th boundary. We define 
\[
	\tDer(\mathbb{K}\pi) = \{(u, u_1,\dotsc,u_n)\in\Der_\mathbb{K}(\mathbb{K}\pi)\times (\mathbb{K}\pi)^n : u(\gamma_j) = [\gamma_j,u_j]\}.
\]
\end{definition}

\begin{definition}
Denote the homotopy class of a simple loop based at $*_j$ representing $\partial_j\Sigma$ abusively by $\partial_j\Sigma$. We define $\Der_\partial(\mathbb{K}\mathscr{G})$ to be the space of derivations $f$ on a $\mathbb{K}$-linear category $\mathbb{K}\mathscr{G}$ with $f(\partial_j\Sigma) = 0$ for all $j$.
\end{definition}

For the definition of a derivation on a $\mathbb{K}$-linear category, see Section 4 of \cite{toyo2}, for example.

\begin{proposition}\label{prop:tder}
We have an identification $\tDer(\mathbb{K}\pi) \cong \Der_\partial(\mathbb{K}\mathscr{G})$.
\end{proposition}
\noindent Proof. Suppose we are given a tangential derivation $(u,u_1,\dotsc,u_n)$. We construct a $\mathbb{K}$-linear category derivation $f$. For a loop $\alpha$ based at $*_0$, we set $f(\alpha) = u(\alpha)$. For a path $\delta$ from $*_j$ to $*_0$ such that $\delta\gamma_j\delta^{-1}$ is the $j$-th boundary loop $\partial_j\Sigma$, we set $\delta^{-1}f(\delta) = u_j$. This is well-defined: for another such path $\delta'$, the loop $\delta'^{-1}\delta$ is based at $*_0$ and we have $(\delta'^{-1}\delta)\gamma_j = \gamma_j(\delta'^{-1}\delta)$. Since $\mathscr{G}$ is a free groupoid, we have $\delta'^{-1}\delta = \gamma_j^{-m}$ for some $m\in \mathbb{Z}$. We compute
\begin{align*}
	\delta'^{-1}f(\delta') &= \gamma_j^{-m}\delta^{-1}f(\delta\gamma_j^m)\\
	&= \gamma_j^{-m}\delta^{-1}f(\delta)\gamma_j^{m} + \gamma_j^{-m} f(\gamma_j^{m})\\
	&= \gamma_j^{-m} u_j \gamma_j^{m} + \sum_{1\leq i\leq m}\gamma_j^{-m+i-1} f(\gamma_j)\gamma_j^{m-i}\\
	&= \gamma_j^{-m} u_j \gamma_j^{m} + \sum_{1\leq i\leq m}\gamma_j^{-m+i-1} (\gamma_ju_j - u_j\gamma_j)\gamma_j^{m-i}\\
	&= \gamma_j^{-m} u_j \gamma_j^{m} - \gamma_j^{-m} u_j \gamma_j^{m} + u_j\\
	&= u_j.
\end{align*}
This shows the well-definedness. We also have
\begin{align*}
	f(\partial_j\Sigma) &= f(\delta\gamma_j\delta^{-1})\\
	&= f(\delta)\gamma_j\delta^{-1} + \delta (\gamma_ju_j - u_j\gamma_j)\delta^{-1} - \delta\gamma_j \delta^{-1}f(\delta)\delta^{-1}\\
	&= \delta u_j \gamma_j\delta^{-1}+ \delta (\gamma_ju_j - u_j\gamma_j)\delta^{-1}- \delta\gamma_j u_j\delta^{-1}\\
	&= 0.
\end{align*}

In the opposite direction, given a derivation $f$ on $\mathbb{K}\mathscr{G}$ with $f(\partial_j\Sigma) = 0$, it automatically yields a derivation $u$ by restricting $f$ to $\mathbb{K}\pi = \mathbb{K}\pi_1(\Sigma,*_0)\subset \mathbb{K}\mathscr{G}$. We set $u_j = \delta^{-1}f(\delta)$. This yields a bijection.\qed

\begin{definition}
Let $\hat{\mathscr{T}} = \hat T(H)\langle\delta_j,\delta_j^{-1}\rangle_{1\leq j\leq n}$ be the $\mathbb{K}$-linear category obtained by formally adjoining the morphisms $\delta_j$ and $\delta_j^{-1}$ from $*_j$ to $*_0$ to the completed algebra $\hat T(H)$ which is the space of endomorphisms at $*_0$. A morphism of $\mathbb{K}$-linear categories
\[
	\theta_\mathrm{exp} \colon \mathbb{K}\mathscr{G} \to \hat{\mathscr{T}},
\]
is given as the unique extension with $\theta_\mathrm{exp}(\delta_j) = \delta_j$. The element $\delta_j$ is defined to have degree 0.
A filtration on $\mathbb{K}\mathscr{G}$ is given by the pull-back by $\theta_\mathrm{exp}$. \\
\end{definition}

\subsection{Connections on Modules}

Next, we define a left $\mathbb{K}\mathscr{G}$-module $\mathscr{N}$ and a connection on $\mathscr{N}$, which eventually give the divergence maps in the last section. Recall that, in \cite{toyo}, the single divergence is based on a left module over a Hopf algebra $A$, from which we obtain a bimodule by a functor $\Phi_A$. In the definition, we used $\mathbb{K}$ as a ground field and also as a left  $A$-module via the augmentation $\varepsilon\colon A\to \mathbb{K}$. In the setting of linear categories, we have a notion of Hopf groupoids, and these two $\mathbb{K}$'s become different spaces, as we shall see below. We freely use the correspondence on modules between a linear category and its category algebra.

\begin{definition}[\cite{fresse}, 9.0.2]
A \textit{Hopf groupoid} over $\mathbb{K}$ is a small category $\mathscr{A}$ enriched over the monoidal category of counital coassociative $\mathbb{K}$-coalgebras together with the anti-homomorphism $\varsigma\colon\mathscr{A} \to \mathscr{A}\op$ of $\mathbb{K}$-linear categories called the antipode, satisfying the Hopf relations.
\end{definition}

\begin{definition}
Let $\mathscr{A}$ be a $\mathbb{K}$-Hopf groupoid. A \textit{Hopf derivation} of $\mathscr{A}$ is a derivation $f\colon \mathscr{A}\to \mathscr{A}$ in the sense of Definition 4.1 in \cite{toyo2}, which is also a coderivation on each Hom-space. We denote by $\Der_\mathrm{Hopf}(\mathscr{A})$ the space of Hopf derivations.
\end{definition}

\begin{definition}\label{def:grpdfilt} 
The $\mathbb{K}$-Hopf groupoid structure of the groupoid algebra $\mathbb{K}\mathscr{G}$ is defined by the following: a counital coassociative (cocommutative) coalgebra structure on the Hom-spaces given by, for $v,w \in \mathrm{Ob}(\mathscr{G}) = V$ and $x\in \mathscr{G}(v,w)$,
\[
	\Delta\colon \mathbb{K}\mathscr{G}(v,w) \to \mathbb{K}\mathscr{G}(v,w)\otimes \mathbb{K}\mathscr{G}(v,w): x \mapsto x\otimes x
\]
and 
\[
	\varepsilon \colon \mathbb{K}\mathscr{G}(v,w) \to \mathbb{K}: x\mapsto 1,
\]
and the antipode is given by
\[
	 \varsigma\colon\mathbb{K}\mathscr{G}(v,w) \to \mathbb{K}\mathscr{G}(w,v): x\to x^{-1}.
\]
The completion $\widehat{\mathbb{K}\mathscr{G}}$ with respect to the filtration in Definition \ref{def:grpdfilt} and its associated graded $\hat{\mathscr{T}}$ are automatically (complete) Hopf groupoids. We set 
\begin{align*}
	\Der_{\mathrm{Hopf},\partial}(\widehat{\mathbb{K}\mathscr{G}}) &= \{f\in\Der_\mathrm{Hopf}(\widehat{\mathbb{K}\mathscr{G}}): f(\partial_j\Sigma) = 0\} \textrm{ and }\\
	\Der^+_{\mathrm{Hopf},\partial}(\widehat{\mathbb{K}\mathscr{G}}) &= \{f\in\Der_{\mathrm{Hopf},\partial}(\widehat{\mathbb{K}\mathscr{G}}): \mathrm{gr}(f) = \id\}.
\end{align*}
\end{definition}

\begin{remark}
The completed free Lie algebra $\hat L(H)$ in the last section is the primitive part of $\hat T(H)$, and $\Der(\hat L(H))$ is canonically isomorphic to the space of (continuous) Hopf derivations $\Der_\mathrm{Hopf}(\hat T(H))$. Since the notion of the primitive part behaves poorly in the setting of multi-objects, we stick to the perspective of Hopf algebras: we use $\Der_{\mathrm{Hopf}}(\widehat{\mathbb{K}\mathscr{G}})$ as a multi-object analogue of $\Der(\hat L(H))$.
\end{remark}

\begin{definition} Let $\mathscr{A}$ be a $\mathbb{K}$-Hopf groupoid.
\begin{itemize}
	\item We define the $\mathbb{K}$-linear category $\mathscr{E}$ by $\mathrm{Ob}(\mathscr{E}) = \mathrm{Ob}(\mathscr{A})$ and $\mathscr{E}(v,w) = \mathbb{K}$ for $v,w\in \mathrm{Ob}(\mathscr{A})$ with the composition map
	\[
		\mathscr{E}(u,v) \otimes \mathscr{E}(v,w) = \mathbb{K}\otimes \mathbb{K} \xrightarrow{\cong} \mathbb{K} =  \mathscr{E}(u,w).
	\]
	Denote the unit $1\in \mathbb{K}=\mathscr{E}(v,w)$ by $e_{vw}$. Then $\varepsilon$ above assembles into a morphism of $\mathbb{K}$-linear categories $\varepsilon\colon \mathscr{A} \to \mathscr{E}$ which reads $\varepsilon(x) = e_{vw}$ in this notation.
	\item Let $\mathbb{K}\mathrm{Ob}(\mathscr{A})$ be a left $\mathscr{E}$-module with $\mathbb{K}\mathrm{Ob}(\mathscr{A})(v) = \mathbb{K} v$, the one-dimensional space generated by the symbol $v$. The action of $\mathscr{E}$ is given explicitly by $e_{vw}\cdot w = v$.
	\item We put $s \colon \mathscr{A}\to \mathbb{K}\mathrm{Ob}(\mathscr{A}) \colon x\mapsto \sum_{v\in\mathrm{Ob}(\mathscr{A})} \varepsilon(x)\cdot v$. We call it the \textit{source map} since for $x\in\mathscr{A}(v,w)$ with $\varepsilon(x) = e_{vw}$, we have $s(x) = v$. This is a left module map over $\varepsilon$.
	\item Recall that $V = \mathrm{Ob}(\mathbb{K}\mathscr{G})$. Regarding $\mathbb{K}V$ as a $\mathbb{K}\mathscr{G}$-module via $\varepsilon$, we define another $\mathbb{K}\mathscr{G}$-module by $\mathscr{N} = \Ker (s \colon \mathbb{K}\mathscr{G}\to \mathbb{K}V)$.
\end{itemize}
\end{definition}

\begin{remark}\label{rem:grpdtwo}
The module $\mathscr{N}$ is a multi-object analogue of $\Ker(\mathbb{K}F_n\xrightarrow{\varepsilon} \mathbb{K})$ in  Proposition 3.8 of \cite{toyo}. Note that two roles carried by the augmentation map on $\mathbb{K}F_n$ (namely, the algebra map $\mathbb{K}F_n \to \mathbb{K}$ and the module map $\mathbb{K}F_n \to \mathbb{K}$ over $\varepsilon$) is now separated into the augmentation $\varepsilon$ and the source map $s$.\\[-7pt]
\end{remark}

Recall that the module over an algebra is said to be \textit{dualisable} if it is finitely generated and projective.

\begin{lemma}
$\mathscr{N}$ is dualisable over $\mathbb{K}\mathscr{G}$.
\end{lemma}
\noindent Proof. We have an identification as $\mathbb{K}\mathscr{G}$-modules:
\[	
	\mathscr{N} \cong \bigoplus_{c = \alpha_i,\beta_i,\gamma_j} \mathbb{K}\mathscr{G}(\argdot,*_0) (1_0 - c) \oplus \bigoplus_{1\leq j\leq n}\mathbb{K}\mathscr{G}(\argdot,*_j)(1_j - \delta_j).
\]
This decomposition is given by the free groupoid version of the Fox derivatives. Then, since the $\mathbb{K}\mathscr{G}$-modules $\mathbb{K}\mathscr{G}(\argdot,*_0)$ and $\mathbb{K}\mathscr{G}(\argdot,*_j)$ are projective and finitely generated (in fact, they are generated by a singe element), the claim follows. For analogous statements and their proofs,  see Proposition 3.8 in \cite{toyo} and Lemma 5.6 in \cite{toyo2}.\qed\\[-7pt]

Now we define a connection on $\mathscr{N}$. For the definition of a connection on a module over a linear category, see Definition 4.8 of \cite{toyo2}.

\begin{definition}
Take the free-generating system $\mathcal{C}$ of $\mathscr{G}$ as in Figure 2 of \cite{toyo2}. We define the connection $\nabla'\!_{\mathcal{C},\mathsf{fr}}\colon \mathscr{N} \to \Omega^1\mathbb{K}\mathscr{G}\otimes_{\mathbb{K}\mathscr{G}}\mathscr{N}$ by 
\begin{align*}
	&\nabla'\!_{\mathcal{C},\mathsf{fr}}(1_0 - \alpha_i) = a_i (d\beta_i)\beta_i^{-1}\otimes (1_0 - \alpha_i),\\
	&\nabla'\!_{\mathcal{C},\mathsf{fr}}(1_0 - \beta_i) = -b_i (d\alpha_i)\alpha_i^{-1} \otimes (1_0 - \beta_i),\\
	&\nabla'\!_{\mathcal{C},\mathsf{fr}}(1_0 - \gamma_j) = 0, \textrm{ and}\\
	&\nabla'\!_{\mathcal{C},\mathsf{fr}}(1_j - \delta_j) = c_j (d\delta_j)\delta_j^{-1}\otimes(1_j - \delta_j).
\end{align*}
Here, $a_i$, $b_i$ and $c_j$ are the rotation numbers of the generators $\alpha_i,\beta_i,\gamma_j\in\mathcal{C}$. Extend $\nabla'\!_{\mathcal{C},\mathsf{fr}}$ to the completion $\hat{\mathscr{N}}$ by continuity, which is a $\widehat{\mathbb{K}\mathscr{G}}$-module. A derivation action of $f\in \Der_\mathrm{Hopf}(\widehat{\mathbb{K}\mathscr{G}})$ on the module $\hat{\mathscr{N}}$ is given by the restriction $f|_{\Ker(s)}$.
\end{definition}

\begin{definition}\
Let $\hat{\mathscr{T}}$ as in Definition \ref{def:grpdfilt} and $\hat{\mathscr{N}}\gr = \Ker(s\colon \hat{\mathscr{T}} \to \mathbb{K}V)$, which is also dualisable over $\hat{\mathscr{T}}$ since it is isomorphic to $\mathscr{N}$ via $\theta_\mathrm{exp}$. We define the connection $\nabla'\!\!_{H,\mathsf{fr}}$ on $\hat{\mathscr{N}}\gr$ by 
\[
	\nabla'\!\!_{H,\mathsf{fr}}(x_i) = \nabla'\!\!_{H,\mathsf{fr}}(y_i) = \nabla'\!\!_{H,\mathsf{fr}}(z_j) = 0\; \textrm{ and } \; \nabla'\!\!_{H,\mathsf{fr}}(1_j - \delta_j) = c_j (d\delta_j)\delta_j^{-1}\otimes(1_j - \delta_j).
\]
This is an extension of the flat connection $\nabla'\!\!_z$ defined in Section 6 of \cite{toyo}. Also this is the associated graded of $\nabla'\!_{\mathcal{C},\mathsf{fr}}$ above since $(d\beta_i)\beta_i^{-1}$ and $(d\alpha_i)\alpha_i^{-1}$ in the definition of $\nabla'\!_{\mathcal{C},\mathsf{fr}}$ have the weight $\geq 1$.\\[-7pt]
\end{definition}

\begin{lemma}\label{lem:sdiv}
We have $\sdiv^\mathsf{fr} = -\Div^{\nabla'\!_{\mathcal{C},\mathsf{fr}}}$ as maps $\Der_{\mathrm{Hopf},\partial}(\widehat{\mathbb{K}\mathscr{G}}) \to |\widehat{\mathbb{K}\mathscr{G}}|$, and $\sdiv^\mathsf{fr}\gr = -\Div^{\nabla'\!\!_{H,\mathsf{fr}}}$ as maps $\Der_{\mathrm{Hopf},\partial}( \hat{\mathscr{T}}) \to | \hat{\mathscr{T}}|$.
\end{lemma}

Before the proof, we look at the counterpart of the connection $\nabla'\!_{\mathcal{C},\mathsf{fr}}$ on a $\mathbb{K}\mathscr{G}$-bimodule $\Omega^1\mathbb{K}\mathscr{G}$. For a $\mathbb{K}$-algebra $A$, we denote the copy of an element $a\in A$ in $A\op$ by $\bar a$. We denote by $A^\mathrm{e}$ the enveloping algebra of $A$.

\begin{definition}
We define the connection $\nabla\!_{\mathcal{C},\mathsf{fr}}\colon\Omega^1\mathbb{K}\mathscr{G}\to \Omega^1\mathbb{K}\mathscr{G}^\mathrm{e} \otimes_{\mathbb{K}\mathscr{G}^\mathrm{e}} \Omega^1\mathbb{K}\mathscr{G} $ on the $\mathbb{K}\mathscr{G}$-bimodule $\Omega^1\mathbb{K}\mathscr{G}$ by
\begin{align*}
	&\nabla\!_{\mathcal{C},\mathsf{fr}}\!\left((d\alpha_i)\alpha_i^{-1} \right) = a_i [(d\beta_i)\beta_i^{-1} - \overline{(d\beta_i)\beta_i^{-1}}]\otimes (d\alpha_i)\alpha_i^{-1},\\
	&\nabla\!_{\mathcal{C},\mathsf{fr}}\!\left((d\beta_i)\beta_i^{-1} \right) = -b_i [(d\alpha_i)\alpha_i^{-1} - \overline{(d\alpha_i)\alpha_i^{-1}}]\otimes (d\beta_i)\beta_i^{-1}\\
	&\nabla\!_{\mathcal{C},\mathsf{fr}}\!\left((d\gamma_j)\gamma_j^{-1} \right) = 0, \textrm{ and}\\
	&\nabla\!_{\mathcal{C},\mathsf{fr}}\!\left((d\delta_j)\delta_j^{-1} \right) = c_j [(d\delta_j)\delta_j^{-1} - \overline{(d\delta_j)\delta_j^{-1}}] \otimes (d\delta_j)\delta_j^{-1}.\\[-7pt]
\end{align*}
\end{definition}

\begin{theorem}\label{thm:cobfr}
Let $\sigma\colon |\mathbb{K}\mathscr{G}| \to \Der_\partial(\mathbb{K}\mathscr{G})$ be the groupoid version of the Kawazumi--Kuno action. For any framing $\mathsf{fr}$, the framed version of the Turaev cobracket $\delta^\mathsf{fr}$ is equal to the composition $-\Div^{\nabla\!_{\mathcal{C},\mathsf{fr}}}\circ\,\sigma$.
\end{theorem}
\noindent Proof. Since the natural map $|\mathbb{K}\mathscr{G}|\to |\widehat{\mathbb{K}\mathscr{G}}|$ is injective, all the computations can be done in the completion, which is further identified with $|\hat{\mathscr{T}}|$ by $\theta_\mathrm{exp}$. For $f\in\Der_\partial(\mathbb{K}\mathscr{G})$, the associated divergence is given by
\begin{align*}
	-\Div^{\nabla\!_{\mathcal{C},\mathsf{fr}}}(f) &= -\Tr(i_f\nabla\!_{\mathcal{C},\mathsf{fr}} - L_f)\\
	&= \sum_{1\leq i\leq g} -|a_i (f(\beta_i)\beta_i^{-1} - \overline{f(\beta_i)\beta_i^{-1}})| + |b_i(f(\alpha_i)\alpha_i^{-1} - \overline{f(\alpha_i)\alpha_i^{-1}})|\\
	&\qquad - \sum_{1\leq j\leq n} |c_j (f(\delta_j)\delta_j^{-1} - \overline{f(\delta_j)\delta_j^{-1}})| + \sum_{c\in\mathcal{C}} |\partial_c(f(c)) - 1\otimes c^{-1}f(c)|\\
	&= \sum_{1\leq i\leq g} a_i|1\wedge f(\beta_i)\beta_i^{-1}| - b_i |1\wedge f(\alpha_i)\alpha_i^{-1}|\\
	&\qquad + \sum_{1\leq j\leq n} c_j|1\wedge f(\delta_j)\delta_j^{-1}| + \sum_{c\in\mathcal{C}} |\partial_c(f(c)) - 1\otimes c^{-1}f(c)|\\
	&= 1\wedge f(\mathbf{p}^\mathsf{fr}) + 1\wedge \mathsf{b}^\mathsf{fr}(f) + \sum_{c=\alpha_i,\beta_i,\gamma_j}|\partial_c(f(c)) - 1\otimes c^{-1}f(c)| + \sum_{1\leq j\leq n} |\partial_{\delta_j}(f(\delta_j)) - 1\otimes \delta_j^{-1}f(\delta_j)|.
\end{align*}
The last term is equal to zero since we have
\[
	\partial_{\delta_j}(f(\delta_j)) - 1\otimes \delta_j^{-1}f(\delta_j) = \partial_{\delta_j}(\delta_ju_j) - 1\otimes u_j = 1\otimes u_j - 1\otimes u_j = 0.
\]
Therefore we have $-\Div^{\nabla\!_{\mathcal{C},\mathsf{fr}}} = \Div^\mathsf{fr}$ on $\Der_\partial(\mathbb{K}\mathscr{G})$, where $\Div^\mathsf{fr}$ is defined in Definition 5.8 of \cite{akkn}. By Theorem 5.16 of \cite{akkn}, we have $\delta^\mathsf{fr} = \Div^\mathsf{fr}\circ\,\sigma$. This completes the proof.\qed\\[-7pt]

The above wording ``counterpart''  is justified by the following.

\begin{lemma}\label{lem:indconn}
The connection $\nabla\!_{\mathcal{C},\mathsf{fr}}$ is induced by $\nabla'\!_{\mathcal{C},\mathsf{fr}}$ in the sense of Definition-Lemma 6.2 in \cite{toyo}.
\end{lemma}
\noindent Proof. Let $A$ be the category algebra of $\mathbb{K}\mathscr{G}$, and $S\cong \prod_V \mathbb{K}$ the corresponding subalgebra generated by the identity morphisms in $\mathbb{K}\mathscr{G}$. The induced connection is only defined in the absolute case $S = \mathbb{K}$ in \cite{toyo}, so we describe the construction for the general case here. The reader may want to refer to Section 6 there beforehand.

First, the source map corresponds to the left $A$-module map $s\colon A\to \mathbb{K}V$, where the corresponding $A$-module to $\mathbb{K}V$ is denoted again by $\mathbb{K}V$. In addition, $\mathscr{N}$ corresponds to $N = \Ker(s\colon A\to \mathbb{K}V)$. Next, consider the functor
\[
	\Phi_{(A,S)} \colon A\textrm{-}\textbf{Mod}\to A^\mathrm{e}\textrm{-}\textbf{Mod}: M \mapsto \Phi_{(A,S)}(M),\;\; \psi\mapsto\psi\otimes\id_{A}\,.
\]
Here we set $\Phi_{(A, S)}(M) = M\otimes_S A$ as $\mathbb{K}$-vector spaces with $M$ regarded as a right $S$-module since $S$ is commutative, and the bimodule structure is given by, for $a,x,y\in A$ and $m\in M$,
\[
	x\cdot(m\otimes a)\cdot y = x^{(1)}m\otimes x^{(2)}ay,
\]
using the coproduct $\Delta(x) = x^{(1)}\otimes x^{(2)}$. This action is well-defined: for two representatives $1_vm\otimes a$ and $m\otimes 1_va$ of the same element in $\Phi_{(A,S)}(M)$, the action of $x\in\mathscr{G}(-,w)$ is computed as, if $w=v$,
\[
	x\cdot (1_vm\otimes a) = xm\otimes xa = x\cdot(m\otimes 1_va),
\]
and if $w\neq v$,
\[
	x\cdot (1_vm\otimes a) = 0 = x\cdot(m\otimes 1_va).
\]
Then, we have isomorphisms of $A$-bimodules: for $x,y\in\mathscr{G}$ and $v\in V$,
\begin{align}\label{eq:hopfisom}
	\Phi_{(A,S)}(A) \cong A\otimes_SA&:  x\otimes y \mapsto x\otimes x^{-1}y\;\textrm{ and}\\
	\Phi_{(A,S)}(\mathbb{K}V) \cong A\hspace{25.5pt}&:  v\otimes x  \mapsto 1_vx\,.\nonumber
\end{align}
The functor $\Phi_{(A, S)}$ is exact since $S$ being separable implies every (left) module over $S$ is projective. By applying it to the short exact sequence
\[
	0 \to N \to A\xrightarrow{s} \mathbb{K}V\to 0,
\]
combined with the two isomorphisms above, we get
\[
	0 \to \Phi_{(A,S)}(N) \to A\otimes_SA \xrightarrow{\mathrm{mult}} A \to 0,
\]
so we have $\Phi_{(A,S)}(N) \cong \Ker(A\otimes_SA \xrightarrow{\mathrm{mult}} A) = \Omega^1_SA$. Furthermore, we have a module map
\[
	j_N\colon N \to \Phi_{(A,S)}(N): n\mapsto n\otimes 1
\]
over the twisted coproduct $\tilde\Delta = (\id\otimes \,\varsigma)\circ \Delta$ by Definition-Lemma 6.1 of \cite{toyo}. The composition of $j_N$ with the isomorphism \eqref{eq:hopfisom} sends $(1_0-c)$ to $(dc)c^{-1}$ for $c\in\mathcal{C}$, and $(1_j - \delta_j)$ to $(d\delta_j)\delta_j^{-1}$, respectively. Since the induced connection $\Phi_{(A,S)}(\nabla'\!_{\mathcal{C},\mathsf{fr}})$ is the natural extension of $\nabla'\!_{\mathcal{C},\mathsf{fr}}$ along $j_N$, we have $\Phi_{(A,S)}(\nabla'\!_{\mathcal{C},\mathsf{fr}}) = \nabla\!_{\mathcal{C},\mathsf{fr}}$ under the identification \eqref{eq:hopfisom}.\qed\\

\noindent \textbf{Proof of Lemma \ref{lem:sdiv}.} Since $\Phi_{(A,S)}(\nabla'\!_{\mathcal{C},\mathsf{fr}}) = \nabla\!_{\mathcal{C},\mathsf{fr}}$, we can apply the Hopf groupoid version of Theorem 6.4 in \cite{toyo} to conclude $\Div^{\nabla\!_{\mathcal{C},\mathsf{fr}}} = |\tilde\Delta|\circ\Div^{\nabla'\!_{\mathcal{C},\mathsf{fr}}}$ on $\Der_\mathrm{Hopf}(\widehat{\mathbb{K}\mathscr{G}})$, where $ |\tilde\Delta|\colon |\widehat{\mathbb{K}\mathscr{G}}| \to |\widehat{\mathbb{K}\mathscr{G}}|\otimes  |\widehat{\mathbb{K}\mathscr{G}\op}|$ is the map induced by the twisted coproduct. Using the augmentation map $|\varepsilon\circ\varsigma|\colon |\widehat{\mathbb{K}\mathscr{G}\op}| \to |\mathscr{E}| \cong \mathbb{K}$, we obtain
\[
	(\id\otimes|\varepsilon\circ\varsigma|)\circ \Div^{\nabla\!_{\mathcal{C},\mathsf{fr}}} = \Div^{\nabla'\!_{\mathcal{C},\mathsf{fr}}}.
\]
We also have $-\Div^{\nabla\!_{\mathcal{C},\mathsf{fr}}} = \Div^\mathsf{fr}$ on $\Der_{\mathrm{Hopf},\partial}(\widehat{\mathbb{K}\mathscr{G}})$ by the proof of Theorem \ref{thm:cobfr}, which is further equal to $|\tilde\Delta|\circ\sdiv^\mathsf{fr}$ by Proposition 4.42 of \cite{akkn}. Hence, we have
\[
	 -\Div^{\nabla'\!_{\mathcal{C},\mathsf{fr}}} = -(\id\otimes|\varepsilon\circ\varsigma|)\circ \Div^{\nabla\!_{\mathcal{C},\mathsf{fr}}} = (\id\otimes|\varepsilon\circ\varsigma|)\circ |\tilde\Delta|\circ\sdiv^\mathsf{fr} = \sdiv^\mathsf{fr}.
\]
For $\nabla'\!\!_{H,\mathsf{fr}}$, it can be similarly done.\qed\\

\subsection{Integration and Connections}

Our goal, for the time being, is to explicitly integrate divergences in terms of connections (and this is essentially done in Proposition 4.39 of \cite{akkn}). First, we define the adjoint action.

\begin{definition}
Let $B$ be a $\mathbb{K}$-algebra and $R$ a subalgebra of $B$.
\begin{itemize}	
	\item An \textit{automorphic action} on a $B$-module $M$ relative to $R$ is a triple $(\Gamma,\varphi,\rho)$ where $\Gamma$ is a group, and
	\[
		\varphi\colon \Gamma \to \Aut_R(B)\;\textrm{ and }\; \rho\colon \Gamma \to \Aut_\mathbb{K}(M)
	\]
	are group homomorphisms satisfying
	\[
		\rho(\gamma)(bm) = \varphi(\gamma)(b)\cdot \rho(\gamma)(m).
	\]
	for $\gamma\in \Gamma$, $b\in B$ and $m\in M$. Here $\Aut_R(B)$ is the set of algebra automorphisms that restricts to the identity on $R$.
	\item Let $\nabla\colon M \to \Omega^1_RB\otimes_ BM$ be an $R$-linear connection on $M$ with an automorphic action. We define the \textit{adjoint action} of $\gamma\in \Gamma$ to $\nabla$ by the composite
\[
	\Ad_\gamma\!\nabla\colon M\xrightarrow{\rho(\gamma^{-1})}M\xrightarrow{\nabla}\Omega^1_RB\otimes_BM\xrightarrow{\varphi(\gamma)\otimes \rho(\gamma)}\Omega^1_RB\otimes_BM\,.
\]
\end{itemize}
\end{definition}

Recall that the space of de Rham forms $\DR_R^\bullet\!B$ is the graded trace space of $\Omega^\bullet_RB$.

\begin{proposition}\label{prop:integral}
Let $B$ and $M$ be equipped with a topology, $\mathfrak d$ a pro-nilpotent Lie algebra and $(\mathfrak d, \varphi, \rho)$ a continuous derivation action $M$. Let $\nabla$ be a trace-flat $R$-linear continuous connection on $M$, namely, $\Tr(\nabla^2) = 0$ in $\DR_R^2\!B$. We denote the integration of $\Div^\nabla\colon \mathfrak d \to |B|$ by $\mathsf{J}^\nabla\colon \Exp(\mathfrak d) \to |B|$. Then, for $G \in \Exp(\mathfrak d)$, we have $d(\mathsf{J}^\nabla(G)) = \Tr((\Ad_G-\id)\nabla)$ in $\DR_R^1\!B$.
\end{proposition}
\noindent Proof. First, the derivation action $(\mathfrak d, \varphi, \rho)$ induces an automorphic action $(\Exp(\mathfrak d), \varphi, \rho)$, and the map $G\mapsto (\Ad_G-\id)\nabla$ on the right-hand side is a group $1$-cocycle. Since the trace map is $\Aut_R(B)$-equivariant, the map $G\mapsto \Tr((\Ad_G-\id)\nabla)$ is also a group $1$-cocycle. On the other hand, the trace-flatness implies $\Div^\nabla$ is a Lie algebra $1$-cocycle:
\[
	\Div^\nabla([f,g]) - f(\Div^\nabla(g)) + g(\Div^\nabla(f)) = i_{\varphi(f)}i_{\varphi(g)}\Tr(\nabla^2)
\]
for $f,g\in\mathfrak d$ where $i$ denotes the contraction. The proof is essentially done in Lemma 4.11 of \cite{toyo}. Therefore $d(\mathsf{J}^\nabla)$ is a group $1$-cocycle. Now, it suffices to show that differentiated Lie algebra $1$-cocycle satisfies the equation
\[
	d(\Div^\nabla(f)) = \Tr\left( \left.\frac{d}{d\varepsilon}\right|_{\varepsilon=0}(\Ad_{e^{\varepsilon f}}-\id) \nabla\right)
\]
for $f\in\mathfrak d$. It is readily seen that $\left.\frac{d}{d\varepsilon}\right|_{\varepsilon=0}(\Ad_{e^{\varepsilon f}}-\id) \nabla$ is equal to the infinitesimal adjoint action
\[
	\ad_f\!\nabla := (L_f\otimes \id + \id\otimes \rho(f))\circ\nabla - \nabla\circ \rho(f),
\]
and that $\Tr(\ad_f\!\nabla) = d(\Div^\nabla(f))$ by Corollary A.5 of \cite{toyo3}. \qed\\

With these preparations, we can finally write down the KV groups in terms of connections. We introduce some more notations.

\begin{definition}\ 
\begin{itemize}
	\item Let $\partial\widehat{\mathbb{K}\mathscr{G}}$ be a topological linear subcategory of $\widehat{\mathbb{K}\mathscr{G}}$ generated by boundary loops $\partial_j\Sigma$, and $\partial\hat{\mathscr{T}}$ its associated graded. We have identifications $|\partial\widehat{\mathbb{K}\mathscr{G}}| = \Big|\sum_j \K{z_j} + \K\xi\Big|$ and  $|\partial\hat{\mathscr{T}}| = \Big|\sum_j \K{z_j} + \K\omega\Big|$.
	\item We set $\Aut^+_{\mathrm{Hopf},\partial}(\hat{\mathscr{T}}) = \{G\in\Aut^+_{\mathrm{Hopf}}(\hat{\mathscr{T}}): G|_{\partial\hat{\mathscr{T}}} = \id \}$ and
	\begin{align*}
		\mathrm{Isom}^+_{\mathrm{Hopf},\partial}(\hat{\mathscr{T}},\widehat{\mathbb{K}\mathscr{G}}) &= \{G\colon\hat{\mathscr{T}}\to \widehat{\mathbb{K}\mathscr{G}} :\textrm{a complete Hopf groupoid isomorphism, }\\
		& \qquad \quad \mathrm{gr}(G) = \id, \textrm{and restricts to } G\colon \partial\hat{\mathscr{T}}\xrightarrow{\cong} \partial\widehat{\mathbb{K}\mathscr{G}}\}.
	\end{align*}
\end{itemize}
\end{definition}

\begin{theorem}\label{thm:kvboundary}
We get the following expression:
\begin{itemize}
	\item $\mathrm{KV}^\mathsf{fr}_{(g,n+1)} = \Big\{G \in\Aut^+_{\mathrm{Hopf},\partial}(\widehat{\mathbb{K}\mathscr{G}}): \Tr(\Ad_G\!\nabla'\!_{\mathcal{C},\mathsf{fr}} - \nabla'\!_{\mathcal{C},\mathsf{fr}}) \in d|\partial\widehat{\mathbb{K}\mathscr{G}}|\Big\}$,
	\item $\mathrm{KRV}^\mathsf{fr}_{(g,n+1)} = \Big\{G \in\Aut^+_{\mathrm{Hopf},\partial}(\hat{\mathscr{T}}):\Tr(\Ad_G\!\nabla'\!\!_{H,\mathsf{fr}} - \nabla'\!\!_{H,\mathsf{fr}}) \in d|\partial\hat{\mathscr{T}}|\Big\}$, and
	\item $\mathrm{SolKV}^\mathsf{fr}_{(g,n+1)} = \Big\{G\in\mathrm{Isom}^+_{\mathrm{Hopf},\partial}(\hat{\mathscr{T}},\widehat{\mathbb{K}\mathscr{G}}): \Tr(\Ad_G\!\nabla'\!\!_{H,\mathsf{fr}} - \nabla'\!_{\mathcal{C},\mathsf{fr}}) \in d|\partial\widehat{\mathbb{K}\mathscr{G}}|\Big\}$.
\end{itemize}
\end{theorem}
\noindent Proof. Since $\log G$ has degree $\geq 1$, traces in the above expressions have no constant terms. Therefore, combining with $\Ker(d\colon |\hat{\mathscr{T}}| \to \DR^1\!\hat{\mathscr{T}}) = \mathbb{K}$, we can work in the space of de Rham $1$-forms without any loss of information. The isomorphisms
\begin{align*}
	\tAut^+(\hat L(H))\cap \{G(\xi)=\xi\} &\cong\Aut^+_{\mathrm{Hopf},\partial}(\widehat{\mathbb{K}\mathscr{G}}),\\
	 \tAut^+(\hat L(H))\cap \{G(\omega)=\omega\} & \cong\Aut^+_{\mathrm{Hopf},\partial}(\hat{\mathscr{T}}),\textrm{ and}\\
	\tAut^+(\hat L(H))\cap \{G(\omega)=\xi\} &\cong \mathrm{Isom}^+_{\mathrm{Hopf},\partial}(\hat{\mathscr{T}},\widehat{\mathbb{K}\mathscr{G}})
\end{align*}
are obtained from Proposition \ref{prop:tder} together with the isomorphism $\theta_\mathrm{exp}$. Next, we can check that two connections $\nabla'\!_{\mathcal{C},\mathsf{fr}}$ and $\nabla'\!\!_{H,\mathsf{fr}}$ are trace-flat. By Lemma \ref{lem:sdiv} and Proposition \ref{prop:integral}, group $1$-cocycles $d(\mathsf{j}^\mathsf{fr}(G))$ and $d(\mathsf{j}^\mathsf{fr}\gr(G))$ are equal to $- \Tr(\Ad_G\!\nabla'\!_{\mathcal{C},\mathsf{fr}} - \nabla'\!_{\mathcal{C},\mathsf{fr}})$ and $-\Tr(\Ad_G\!\nabla'\!\!_{H,\mathsf{fr}} - \nabla'\!\!_{H,\mathsf{fr}})$, respectively. This completes the proof for $\mathrm{KV}^\mathsf{fr}_{(g,n+1)}$ and $\mathrm{KRV}^\mathsf{fr}_{(g,n+1)}$.

For $\mathrm{SolKV}^\mathsf{fr}_{(g,n+1)}$, we first compute 
\begin{align*}
	( \nabla'\!\!_{H,\mathsf{fr}} - \Ad_{\theta_\mathrm{exp}}\!\nabla'\!_{\mathcal{C},\mathsf{fr}})(1_0 - e^{x_i}) &= \nabla'\!\!_{H,\mathsf{fr}} (1_0 - e^{x_i}) - \theta_\mathrm{exp}(\nabla'\!_{\mathcal{C},\mathsf{fr}})(1_0-\alpha_i)\\
	&= d\left(\frac{1_0-e^{x_i}}{x_i}\right)\otimes x_i - a_i (de^{y_i})e^{-y_i}\otimes (1_0 - e^{x_i})\\
	&= \left( d\left(\frac{1_0-e^{x_i}}{x_i}\right)\cdot\frac{x_i}{1_0-e^{x_i}} - a_i (de^{y_i})e^{-y_i}\right)\otimes (1_0 - e^{x_i}),
\end{align*}
whose coefficient is equal to $|dr(x_i) - a_idy_i|$ in $\DR^1\hat{\mathscr{T}}$. It is similar for $y_i$ and $z_j$. For $\delta_j$, we have
\[
	( \nabla'\!\!_{H,\mathsf{fr}} - \Ad_{\theta_\mathrm{exp}}\!\nabla'\!_{\mathcal{C},\mathsf{fr}})(1_j-\delta_j)  = 0
\]
by definition. Therefore, we have
\begin{align*}
	\Tr(\nabla'\!\!_{H,\mathsf{fr}} - \Ad_{\theta_\mathrm{exp}}\!\nabla'\!_{\mathcal{C},\mathsf{fr}}) = \sum_i d|r(x_i) - a_iy_i + r(y_i) + b_ix_i| + \sum_j d|r(z_j)| = d|\mathbf r' - \mathbf{p}^\mathsf{fr}|.
\end{align*}
Now for $G\in\mathrm{Isom}^+_{\mathrm{Hopf},\partial}(\hat{\mathscr{T}},\widehat{\mathbb{K}\mathscr{G}})$, we put $F = \theta_\mathrm{exp}\circ G\in \Aut^+_{\mathrm{Hopf}}(\hat{\mathscr{T}})$. We compute
\begin{align*}
	\Tr(\Ad_G\!\nabla'\!\!_{H,\mathsf{fr}} - \nabla'\!_{\mathcal{C},\mathsf{fr}}) &= \Tr(\Ad_{\theta_\mathrm{exp}^{-1}\circ F}\!\nabla'\!\!_{H,\mathsf{fr}} - \Ad_{\theta_\mathrm{exp}^{-1}}\!\nabla'\!\!_{H,\mathsf{fr}}) + \Tr(\Ad_{\theta_\mathrm{exp}^{-1}}\!\nabla'\!\!_{H,\mathsf{fr}} - \nabla'\!_{\mathcal{C},\mathsf{fr}})\\
	&=  \theta_\mathrm{exp}^{-1}\left[\Tr(\Ad_F\!\nabla'\!\!_{H,\mathsf{fr}} - \nabla'\!\!_{H,\mathsf{fr}}) + \Tr(\nabla'\!\!_{H,\mathsf{fr}} - \Ad_{\theta_\mathrm{exp}}\!\nabla'\!_{\mathcal{C},\mathsf{fr}})\right]\\
	&= \theta_\mathrm{exp}^{-1}\left[d(\mathsf{J}^{\nabla'\!\!_{H,\mathsf{fr}}}(F) )+  d|\mathbf r' - \mathbf{p}^\mathsf{fr}|\right]\\
	&= \theta_\mathrm{exp}^{-1}\left[- d(\mathsf{j}^\mathsf{fr}\gr(F)) +  d|\mathbf r' - \mathbf{p}^\mathsf{fr}|\right].
\end{align*}
Using $\mathbf r - \mathbf r' \in \partial\widehat{\mathbb{K}\mathscr{G}}$, this is exactly the defining equation of the KV associators. \qed\\

\section{Rewriting Rules}\label{sec:rewriting}
Consider the case $\Sigma = \Sigma_{g,0}$. The first homology $H = H_1(\Sigma;\mathbb{K})$ is a symplectic space and we have $\omega = \sum_{1\leq i\leq g}[x_i,y_i]$. We put $T(H)_\omega = T(H)/\omid$, which is the quotient of (non-completed) $T(H)$ by the two-sided ideal $\omid$. The purpose of this section is to construct a basis of $|T(H)_\omega|$ for the later sections. See Theorem \ref{thm:basis} for the final result.\\

A basis is constructed by means of rewriting rules. First of all, consider the following (one-way) rewriting in $|T(H)|$: for any monomial $b$ in $T(H)$,
\[
	|by_gx_g| \xmapsto{\rho} |b(x_gy_g + \omega') |.
\]
Here we put $\omega' =  \sum_{1\leq i <g}[x_i,y_i] $. This rewriting may not terminate, so we introduce the \textit{irregularity} to measure the number of steps for the termination. 

\begin{definition}
Let $C$ be the set of all (non-commutative) monomials over $H$ so that $C$ is a basis of $|T(H)|$. For $w \in C$, its irregularity $\mathrm{irr}(w)$ is defined by the following steps. For each appearance $\eta$ of $y_g$ in $w$, compute $t_\eta\in\mathbb{Z}_{\geq 0}\cup \{+\infty\}$ by the following procedure:
\begin{enumerate}[(1)]
	\item Start from the letter $\eta$ in $w$. Set $t_\eta = 0$;
	\item See the next (cyclically adjacent to the right) letter $h$. If $h=x_g$, add $1$ to $t_\eta$. If $h=y_g$, add $0$ to $t_\eta$. Otherwise, terminate the procedure;
	\item Go back to (2).
\end{enumerate}
If this procedure terminates, we define $t_\eta$ as the final value we obtained. If not, we set $t_\eta$ as the supremum of the values attained during the procedure, either $0$ or $+\infty$. We define the irregularity of $w$ by
\[
	\irr(w) = \sum_\eta t_\eta \in\mathbb{Z}_{\geq 0}\cup \{+\infty\}.
\]
If there is no $y_g$'s in $w$, the sum is empty and $\irr(w) = 0$. For a polynomial $a = \sum_{w\in C} c_w w\in |T(H)|$ with $c_w\in\mathbb{K}$, we set $\irr(a) = \sup_{w:c_w\neq 0}( \irr(w)) \in \mathbb{Z}_{\geq 0}\cup \{+\infty\}$.
\end{definition}

\begin{example} Assume $g\geq 2$.
\begin{itemize}
	\item If $w$ is the empty word (corresponding to $|1|\in |T(H)|$), we have $\irr(w) = 0$. 
	\item Consider the word $w = |y_gx_gx_1|$. There is only one $y_g$ in $w$, so we set $\eta$ as that. Following the procedure above, we get $t_\eta = 1$. We have $\irr(w) = 1$.
	\item Consider the word $w = |y_gx_gx_gy_gx_gy_1|$. There are two $y_g$'s in $w$, and each yields $t_\eta = 3$ and $1$, respectively. We have $\irr(w) = 3 + 1 = 4$.
	\item Consider the word $w = |y_gx_gx_g|$. There is only one $y_g$ in $w$. The above procedure does not terminate, encountering infinitely many $x_g$'s throughout. We have $\irr(w) = +\infty$.
	\item Consider the word $w = |y_g|$. There is only one $y_g$ in $w$. The above procedure does not terminate while encountering no $x_g$'s. We have $\irr(w) = 0$.\\[-7pt]
\end{itemize}
\end{example}

\begin{lemma}\label{lem:irrinf}
For $w\in C$, we have $\irr(w)=+\infty$ if and only if all of the followings hold:
\begin{enumerate}[(1)]
	\item $w$ consists only of $x_g$'s and $y_g$'s;
	\item $w$ contains both $x_g$ and $y_g$.
\end{enumerate}
\end{lemma}
\noindent Proof. Assume (1) and (2). By condition (1), the procedure above does not terminate for any choice of $\eta$. From (2), we can see that $t_\eta = +\infty$.

In the other direction, let us assume the negation of (1). Then, the above procedure terminates, and $t_\eta$ is always finite for any $\eta$. If we assume (1) and the negation of (2), the only possible cyclic word is either $|x_g^m|$ or $|y_g^m|$ for some $m\geq 0$. In any case, we have $\irr(w) = 0<+\infty$.\qed\\[-7pt]

\begin{lemma}
Rewriting by $\rho$ reduces the irregularity: for $a\in |T(H)|$ with $0<\irr(a)<+\infty$, we can obtain $a'\in |T(H)|$ by applying $\rho$ such that $\irr(a') < \irr(a)$.
\end{lemma}
\noindent Proof. It suffices to see the case $a = w\in C$. By the assumption $\irr(w)>0$, there is a sequence $y_gx_g$ in $w$. We put $w = |by_gx_g|$ for some monomial $b\in T(H)$. By the assumption $\irr(w)<+\infty$, $b$ contains $x_i$ or $y_i$ for $1\leq i <g$. By applying $\rho$, we obtain a polynomial $|b(x_gy_g + \omega')|$, whose irregularity is evaluated as
\[
	\irr(|b(x_gy_g + \omega')|) = \max\{\irr(|bx_gy_g|),\irr(|bx_iy_i|)\} = \irr(|bx_gy_g|) =  \irr(|by_gx_g|) - 1.
\]
This completes the proof.\qed\\[-7pt]

\begin{corollary}
For $a\in |T(H)|$, rewriting of the polynomial $a$ by $\rho$ terminates if and only if $\irr(a)<+\infty$.\qed\\[-7pt]
\end{corollary}

\begin{lemma}\label{lem:confluent}
Rewriting by $\rho$ is confluent: if $a'_1$ and $a'_2$ are obtained from $a\in |T(H)|$ by applying $\rho$ several times, there exists $a''\in |T(H)|$ by also applying $\rho$ several times that can be obtained either from $a'_1$ or $a'_2$.
\end{lemma}
\noindent Proof. It suffices to show where $a'_1\neq a'_2$ are obtained by only one application of $\rho$. By the assumption,  there are two ways of applying $\rho$ to $a$: one leads to $a'_1$ and the other to $a'_2$. In this case, we can find distinct sequences of $y_gx_g$ in $a$ (either in the same monomial or different monomials) since it is impossible to share a part of sequence $y_gx_g$ without sharing the entirety of $y_gx_g$. Therefore, we can apply $\rho$ to the other sequence $y_gx_g$ in $a'_1$ and $a'_2$ respectively to obtain a common polynomial $a''$.\qed\\[-7pt]

\begin{corollary}\label{cor:terminates}
If the rewriting terminates, it leads to a unique polynomial. For $a\in |T(H)|$ with $\irr(a)$, we denote the thusly obtained element by $\rho(a)$.\qed\\
\end{corollary}

Let us construct the basis of $|T(H)_\omega|$. From now on, the word ``relation'' refers to an expression of the form $\sum_\lambda c_\lambda w_\lambda\in|\omid|$ where $c_\lambda\in\mathbb{K}$ and $w_\lambda\in C$. First of all, since the two-sided ideal $\omid$ of $T(H)$ is generated by the sole element $\omega$, all the relations amongst elements of $C$ is given as a linear combination of the relation $|b\omega |= |b(x_gy_g + \omega')| - |by_gx_g|\in |\omid|$ for $b\in T(H)$, which come from the rule $\rho$ as a linear combinations of $w - a \in |\omid|$ where $w\in C$ and $a$ is a polynomial obtained by applying $\rho$ to $w$. 

If the rewriting terminates, we have a unique expression by the corollary above, so we have the normal form. In the case of $\irr(w)=+\infty$ for $w\in C$, however, several application of $\rho$ to $w$ can bring back to $w$ itself. For example, applying $\rho$ to the word $w = |y_gx_gx_g|$, we obtain
\[
	a = |(x_gy_g + \omega')x_g| = w+ \sum_{1\leq i <g}|[x_i,y_i]x_g|.
\]
Since $a$ and $w$ give the same element in $|T(H)_\omega|$, this yields a relation $\sum_{1\leq i <g}|[x_i,y_i]x_g| \in|\omid| $, which is a consequence of the rule $\rho$ but does not contain the sequence $y_gx_g$ in its expression. We deal with this type of relations from now on.

We need some preparation.

\begin{definition}
Let $\mathcal{R}$ be a directed graph defined by the following:
\begin{itemize}
	\item Each vertex is labelled by $w\in C$ with $\irr(w) = +\infty$, which is a monomial consisting of $s$ $x_g$'s and $t$ $y_g$'s for some $s,t\geq 1$.
	\item For one application of $\rho$ to $w$ with the result $w' + a$, where $w'\in C$ with $\irr(w) = +\infty$ and $a\in |T(H)|$ with $\irr(a)<+\infty$, we draw one directed edge from $w$ to $w'$. 
\end{itemize}
For an edge $e\colon w\mapsto w'+a$, we define the \textit{holonomy} $\mathrm{hol}(e)$ along the edge (with respect to the fixed gauges $C$) to be $\rho(a)$, which is well-defined by Corollary \ref{cor:terminates}.

For $s,t\geq 1$, $\mathcal{R}^{(s,t)}$ denotes the full subgraph of $\mathcal{R}$ whose vertex labels are comprised of $s$ $x_g$'s and $t$ $y_g$'s. We have a decomposition $\mathcal{R} = \bigsqcup_{s,t\geq 1}\mathcal{R}^{(s,t)}$ as a directed graph.
\end{definition}

A \textit{directed path} is a (possibly empty) sequence of composable directed edges, and a \textit{loop} is a directed path whose starting point coincides with the endpoint, which we call the \textit{base point} of the loop. A \textit{free loop} is an equivalence class of based loops generated by cyclic permutations of edges. The holonomy along a directed path is defined as the sum of the holonomies of edges comprising that path. For a pair of directed paths $w \mapsto w'+ a$ and $w\mapsto w'+ a'$ in $\mathcal{R}$, we obtain a relation $a-a'\in |\omid|$. This is equivalent to the relation $\rho(a) -\rho(a')\in |\omid|$ between their holonomies since we already know $\rho(a) - a \in |\omid|$ for any $a$ with $\irr(a)<+\infty$. All the relations are linearly generated by these, as mentioned above.

\begin{lemma}\label{lem:dirconn}
For two vertices $w,w'$ in the same subgraph $\mathcal{R}^{(s,t)}$, there is a directed path from $w$ to $w'$.
\end{lemma}
\noindent Proof. One application of $\rho$ to $w$ swaps one $y_g$ in $w$ with the $x_g$ immediately after this $y_g$. By applying $\rho$ many times, we can find a directed path from $w$ to $|x_g^sy_g^t|$. Now we can re-distribute $y_g$'s by applying $\rho$ many more times to reach the vertex $w'$.\qed\\[-7pt]

\begin{lemma}
All the relations in $C$ are linearly generated by the relations obtained from loops in $\mathcal{R}$.
\end{lemma}
\noindent Proof. Given a pair of directed paths $w \mapsto w'+ a$ and $w\mapsto w'+ a'$ considered above, we take a directed path $w' \mapsto w + a''$, whose existence is guaranteed by the lemma above. Then, we obtain a pair of composite rewritings
\[
	(w \mapsto w'+a \mapsto w + a + a''\;,\;w \mapsto w'+a' \mapsto w + a' + a'').
\]
This results in the relation $(a+a'') - (a'+a'')\in|\omid|$, which is exactly the relation $a-a'\in |\omid|$ obtained above. 
On the other hand, denoting the constant loop at $w$ by $\mathsf{c}_w$, the pairs of loops
\[
	(w \mapsto w'+a \mapsto w + a + a'', \mathsf{c}_w)\;\textrm{ and }\; (w \mapsto w'+a' \mapsto w + a + a'', \mathsf{c}_w)
\]
yield the relations $a+a'', a'+a''\in|\omid|$, from which the relation $a-a'\in|\omid|$ follows. Therefore, relations coming from loops in $\mathcal{R}$ are enough to deduce all the relations. \qed\\[-7pt]

\begin{lemma}
The holonomy along a loop does not depend on the base point.
\end{lemma}
\noindent Proof. For a loop traversing $w_0,w_1,\dotsc,w_r = w_0$ in this order with base point $w_i$, its holonomy is given by the sum $\sum_{0\leq i<r}a_i$, where $a_i$ is the holonomy along the edge from $w_i$ to $w_{i+1}$. This is clearly independent of the base point $w_i$.\qed\\[-7pt]

\begin{definition}
Let $\tilde{\mathcal{R}}^{(s,t)}$ be the directed graph defined by the following:
\begin{itemize}
	\item Each vertex is labelled by a cyclic sequence $\tilde w$ of letters $x_g^{[i]}$ and $y_g^{[j]}$ for $i\in\mathbb Z/s\mathbb Z$ and $j\in\mathbb Z/t\mathbb Z$ such that they appear exactly once while respecting the cyclic orders $(x_g^{[1]},\dotsc,x_g^{[s]})$ and $(y_g^{[1]},\dotsc,y_g^{[t]})$. The letters $x_g^{[i]}$ and $y_g^{[j]}$ are distinguished copies of $x_g$ and $y_g$, respectively.
	\item For one application of $\rho$ to a cyclic sequence $\tilde w$ with the result $\tilde{w}' + a$, where $w'\in C$ with $\irr(\tilde{w}')=+\infty$ and $a\in |T(H)|$ with $\irr(a)<+\infty$, we draw one edge from $\tilde w$ to $\tilde{w}'$. 
\end{itemize}
The natural quotient map $\tilde{\mathcal{R}}^{(s,t)} \twoheadrightarrow \mathcal{R}^{(s,t)}$ is induced by the re-labelling action of $\mathbb{Z}/s\mathbb{Z}\times \mathbb Z/t\mathbb Z$, which sends $x_g^{[i]}$ and $y_g^{[j]}$ appearing in labels to $x_g$ and $y_g$ respectively.
\end{definition}

For a loop $\ell\colon w \mapsto w + a$ in $\mathcal{R}^{(s,t)}$, its $st$-fold power $\ell^{st}$ lifts to a loop in $\tilde{\mathcal{R}}^{(s,t)}$ and results in a rewriting $\tilde w \mapsto \tilde w + st\cdot a$, which gives a relation $st\cdot a\in|\omid|$. Since we have assumed $\mathrm{char}(\mathbb{K}) = 0$, this is equlvalent to $a\in|\omid|$. For this reason, it is enough to only consider loops in $\tilde{\mathcal{R}}^{(s,t)}$.

\begin{lemma}
Let $w=|(y_g^{[1]}x_g^{[i_1]}\cdots x_g^{[i_2-1]})\cdots (y_g^{[j]}x_g^{[i_j]}\cdots x_g^{[i_{j+1}-1]})\cdots (y_g^{[t]}x_g^{[i_t]}\cdots x_g^{[i_1-1]})|$ be a vertex of $\tilde{\mathcal{R}}^{(s,t)}$. Then, we have a bijection
\begin{align*}
	\{\textrm{edges in }\tilde{\mathcal{R}}^{(s,t)}\textrm{ starting from } w\} \leftrightarrow \{j\in\mathbb Z/t\mathbb Z: i_j<i_{j+1}\}.
\end{align*}
\end{lemma}
\noindent Proof. An edge in $\tilde{\mathcal{R}}^{(s,t)}$ starting from $w$ can be written as the rewriting
\begin{align*}
	m_j&\colon |(y_g^{[1]}x_g^{[i_1]}\cdots x_g^{[i_2-1]})\cdots (y_g^{[j]}x_g^{[i_j]}\cdots x_g^{[i_{j+1}-1]})\cdots (y_g^{[t]}x_g^{[i_t]}\cdots x_g^{[i_1-1]})|\\
	&\quad \mapsto |(y_g^{[1]}x_g^{[i_1]}\cdots x_g^{[i_2-1]})\cdots (x_g^{[i_j]}y_g^{[j]}x_g^{[i_j+1]}\cdots x_g^{[i_{j+1}-1]})\cdots (y_g^{[t]}x_g^{[i_t]}\cdots x_g^{[i_1-1]})| + (\textrm{terms with }\irr<+\infty)\,.
\end{align*}
for some $j\in\mathbb Z/t\mathbb Z$. This moves $x_g^{[i_j]}$ to the left side of $y_g^{[j]}$ without changing other parts, and such a rewriting is possible only when $i_j<i_{j+1}$. The correspondence $m_j\leftrightarrow j$ yields the bijection we want.\qed\\

As a consequence of the above lemma, any directed path in $\tilde{\mathcal{R}}^{(s,t)}$ can be expressed by a sequence $m_{j_1}\cdots m_{j_r}$, read from left to right. At this point, it should be convenient to introduce a diagrammatical explanation for the moves $m_j$. As in Figure \ref{fig:beadscorresp}, consider a circle whose arc is separated in $t$ parts by labelled partitions and $s$ beads are distributed over the segments. Each partition corresponds to one $y_g^{[j]}$, and each bead correspond to one $x_g^{[i]}$. One configuration of beads exactly specifies a unique vertex in $\tilde{\mathcal{R}}^{(s,t)}$. The rewriting $m_j$ removes the $i_j$-th bead, which is immediately right to the $j$-th partition, and adds to immediately left of the partition (Figure \ref{fig:beadsmove}). 

\begin{figure}[t]
\begin{minipage}[b]{0.5\columnwidth}
\centerline{\includegraphics[scale=1]{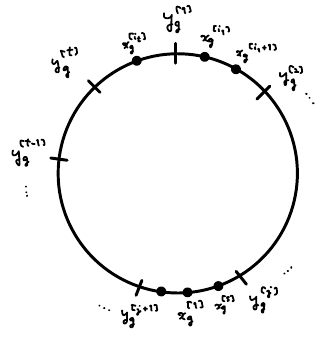}}
\caption{The configuration of beads corresponding to the vertex $|(y_g^{[1]}x_g^{[i_1]}\cdots x_g^{[i_2-1]})\cdots (y_g^{[t]}x_g^{[i_t]}\cdots x_g^{[i_1-1]})|$.}
\label{fig:beadscorresp}
\end{minipage}
\begin{minipage}[b]{0.45\columnwidth}
\centerline{\includegraphics[scale=1]{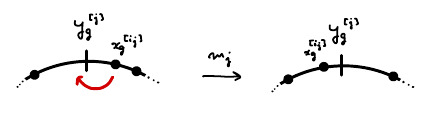}}\vspace{55pt}
\caption{The move $m_j$.\\\ }
\label{fig:beadsmove}
\end{minipage}
\end{figure}

\begin{lemma}
Two directed paths $m_im_j$ and $m_jm_i$ have the same endpoints, and we have $\mathrm{hol}(m_im_j) =\mathrm{hol}(m_jm_i)$.
\end{lemma}
\noindent Proof. Let $\tilde w = |b(y_g^{[i]}x_g^{[p]})b'(y_g^{[j]}x_g^{[q]})|$ for some $p,q$ and $b,b'\in T(H)$. Then, the rewriting of $w'$ along $m_im_j$ reads
\begin{align*}
	|b(y_g^{[i]}x_g^{[p]})b'(y_g^{[j]}x_g^{[q]})| &\xmapsto{m_i} |b(x_g^{[p]}y_g^{[i]})b'(y_g^{[j]}x_g^{[q]})| + |b\omega'b'(y_g^{[j]}x_g^{[q]})|  \\
	&\xmapsto{m_j} |b(x_g^{[p]}y_g^{[i]})b'(x_g^{[q]}y_g^{[j]})| + |b(x_g^{[p]}y_g^{[i]})b'\omega'| + |b\omega'b'(y_g^{[j]}x_g^{[q]})|,
\end{align*}
whose holonomy is 
\[
	\rho\big(|b(x_gy_g)b'\omega'| + |b\omega'b'(y_gx_g)|\big) = \rho\big(|b(x_gy_g)b'\omega'| + |b\omega'b'(x_gy_g + \omega')|\big)\,.
\]
On the other hand, the rewriting of $w'$ along $m_jm_i$ reads
\begin{align*}
	|b(y_g^{[i]}x_g^{[p]})b'(y_g^{[j]}x_g^{[q]})| &\xmapsto{m_j} |b(y_g^{[i]}x_g^{[p]})b'(x_g^{[q]}y_g^{[j]})| +  |b(y_g^{[i]}x_g^{[p]})b'\omega')|\\
	&\xmapsto{m_i} |b(x_g^{[p]}y_g^{[i]})b'(x_g^{[q]}y_g^{[j]})| + |b\omega'b'(x_g^{[q]}y_g^{[j]})| +  |b(y_g^{[i]}x_g^{[p]})b'\omega'|
\end{align*}
whose holonomy is 
\[
	\rho\big(|b\omega'b'(x_gy_g)| +  |b(y_gx_g)b'\omega'|\big) = \rho\big(|b\omega'b'(x_gy_g)| +  |b(x_gy_g + \omega')b'\omega'|\big) \,.
\]
The endpoints and their holonomies are equal, as is claimed.\qed\\

From now on, we study the set of directed paths in modulo the relation $m_im_j= m_jm_i$.

\begin{definition}
For a vertex $w$ of $\tilde{\mathcal{R}}^{(s,t)}$, we set $P^{(s,t)}_w = \{\textrm{directed paths in }\tilde{\mathcal{R}}^{(s,t)} \textrm{ starting from } w\}$. We say a directed path $p_1$ is a \textit{deformation} of $p_2$ if $p_1$ and $p_2$ give the same element in $P^{(s,t)}_w/\langle m_im_j=m_jm_i\rangle$, which is the quotient by the relations multiplicatively generated by $m_im_j=m_jm_i$ for all $i,j$, whenever defined.
\end{definition}
\noindent By the lemma above, the holonomy map is well-defined on the set $P^{(s,t)}_w/\langle m_im_j=m_jm_i\rangle$.

\begin{definition}
We define the \textit{counting number map} by 
\begin{align*}	
	\mathbf{c}\colon P^{(s,t)}_w\to (\mathbb{Z}_{\geq 0})^t\colon p \mapsto (\mathbf{c}(p,j))_{j\in \mathbb{Z}/t\mathbb{Z}}
\end{align*}
where $\mathbf{c}(p,j) = (\textrm{the number of beads jumped over }y_g^{[j]} \textrm{ along }p) = (\textrm{the number of }m_j\textrm{'s in }p)$. \end{definition}

\begin{lemma}\label{lem:rot}
The counting number map induces an injection $P^{(s,t)}_w/\langle m_im_j=m_jm_i\rangle \to (\mathbb{Z}_{\geq 0})^t$.
\end{lemma}
\noindent Proof. Let $w=|(y_g^{[1]}x_g^{[i_1]}\cdots x_g^{[i_2-1]})\cdots (y_g^{[t]}x_g^{[i_t]}\cdots x_g^{[i_1-1]})|$. Suppose that two directed paths $p$ and $p'$ have the same counting number: $\mathbf{c}(p) = \mathbf{c}(p')$. Put $\nu = \sum_{1\leq j\leq t} \mathbf{c}(p,j)$, which is the total number of $m_j$'s appearing in each of $p$ and $p'$. We can write as $p = m_{j_1}\cdots m_{j_\nu}$ and $p' = m_{j'_1}\cdots m_{j'_\nu}$. By the assumption $\mathbf{c}(p) = \mathbf{c}(p')$, $p$ and $p'$ has the same number of $m_j$'s for each $j$. 

We prove by induction on $\nu$. If $\nu=0$, they are constant loops at $w$, so they are equal. Now suppose $\nu \geq 1$. Let $\alpha \in \{1,\dotsc, \nu\}$ such that $m_{j'_\alpha}$ is the leftmost appearence of $m_{j_1}$ in $p'$. Since $m_{j_1}$ is the first edge of $p$, we have $i_{j_1}<i_{j_1+1}$, which says that there is at least one bead between $j_1$-th and $(j_1+1)$-th partitions. For $1\leq \beta<\alpha$, we have $j'_\beta\neq j_1$ so the move $m_{j'_\beta}$ does not reduce the number of beads between $j_1$-th and $(j_1+1)$-th partitions. Therefore, we have $m_{j'_\beta}m_{j'_\alpha} = m_{j'_\alpha}m_{j'_\beta}$ for all $1\leq \beta<\alpha$; hence we can assume $j'_1 = j_1$. Now we can use the induction hypothesis to $m_{j_2}\cdots m_{j_\nu}$ and $m_{j'_2}\cdots m_{j'_\nu}$ to conclude that they are equal modulo the relation $m_im_j=m_jm_i$. This shows $p = p'$ in $P_w^{(s,t)}$.\qed\\

We have a preferred base point $*_{s,t} = |y_g^{[1]}\cdots y_g^{[t]}x_g^{[1]}\cdots x_g^{[s]}|$ in $\tilde{\mathcal{R}}^{(s,t)}$ together with a based loop $m_t^sm_{t-1}^s\cdots m_1^s$, which first passes the edge specified by $t\in\mathbb Z/t\mathbb Z$ exactly $s$ times, then $(t-1)\in\mathbb Z/t\mathbb Z$ exactly $s$ times, and so on, starting at the above base point.

\begin{figure}[t]
\centerline{\includegraphics[scale=1]{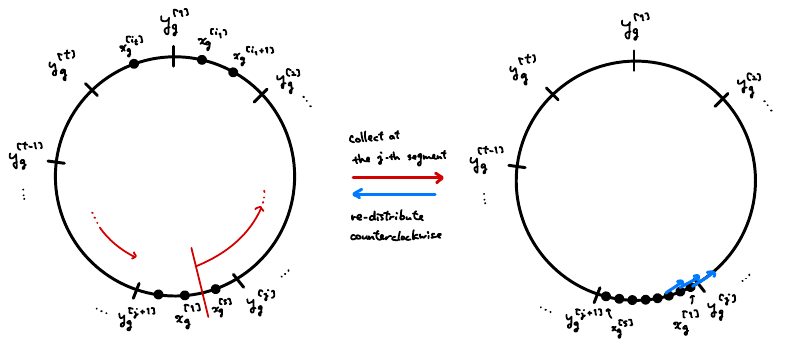}}
\caption{The collection move and the re-distribution move.}
\label{fig:beadsdist}
\end{figure}

\begin{lemma}\label{lem:deform}
For any non-constant loop $\ell$ in $\tilde{\mathcal{R}}^{(s,t)}$ represented by a composition of $m_j$'s, there is a sequence of deformations and cyclic permutations of $m_j$'s which converts $\ell$ to $(m_t^sm_{t-1}^s\cdots m_1^s)^n$ for some $n\geq 1$.
\end{lemma}
\noindent Proof. Let $w=|(y_g^{[1]}x_g^{[i_1]}\cdots x_g^{[i_2-1]})\cdots (y_g^{[t]}x_g^{[i_t]}\cdots x_g^{[i_1-1]})|$ be a base point of $\ell$. Since this is non-constant, at least one bead is moved by $\ell$. However, since $\ell$ is a loop, that one bead has to go around the circle to get back to the original position; this forces every bead, especially $x_g^{[1]}$, to go around the circle. We have $\mathbf{c}(\ell) = (ns,ns,\dotsc, ns)$ for some $n\in\mathbb{Z}_{>0}$. 

Now consider the \textit{collection} move and the \textit{re-distribution} move depicted in Figure \ref{fig:beadsdist}. Let $j$ be such that $i_j \leq 1< i_{j+1}$. The collection move is given by $m_w^\mathrm{col} = m_j^{r_{w,j}}m_{j-1}^{r_{w,j-1}}\cdots m_{j+1}^{r_{w,j+1}}$ where $r_{w,k} = s + 1 - i_k \in \{0,\dotsc,s-1\}$. This move gathers all the beads in $j$-th segment satisfying $i_j \leq 1< i_{j+1}$ by pushing the red line in the figure around the circle so that it ends at the vertex $|y_g^{[j+1]}\cdots y_g^{[j]}x_g^{[1]}\cdots x_g^{[s]}|$. The re-distribution move is given by $m_w^\mathrm{red} = m_j^{s- r_{w,j}}m_{j-1}^{s- r_{w,j-1}}\cdots m_{j+1}^{s- r_{w,j+1}}$. This move sends the beads in bulk along the circle until they reach the original configuration $w$. Then, the composition $\ell' = (m_w^\mathrm{col}m_w^\mathrm{red})^n$ is a loop, and the counting number is $\mathbf{c}(\ell') = (ns, \dotsc, ns) = \mathbf{c}(\ell)$. By Lemma \ref{lem:rot}, $\ell$ is a deformation of $\ell'$.

Put $\ell'' =  (m_w^\mathrm{red} m_w^\mathrm{col})^n$, which is a cyclic permutation of $\ell'$ based at $|y_g^{[j+1]}\cdots y_g^{[j]}x_g^{[1]}\cdots x_g^{[s]}|$. On the other hand, $(m_j^sm_{j-1}^s\cdots m_{j+1}^s)^n$ also has the same base point and the counting number. Again by Lemma \ref{lem:rot}, $\ell''$ is a deformation of $(m_j^sm_{j-1}^s\cdots m_{j+1}^s)^n$. Finally, this is a cyclic permutation of the one in the claim. \qed

\begin{corollary}\label{cor:deform}
For any loop $\ell$ in $\tilde{\mathcal{R}}^{(s,t)}$, we have $\mathrm{hol}(\ell) = n\cdot \mathrm{hol}(m_t^sm_{t-1}^s\cdots m_1^s)$ for some $n\in\mathbb{Z}$.\qed\\[-7pt]
\end{corollary}

\begin{computation}\label{comp:hol}
We compute $h = \mathrm{hol}(m_t^sm_{t-1}^s\cdots m_1^s)$. The path $m_j^s$ projects to the same loop $\ell_{(s,t)}$ in $\mathcal{R}^{(s,t)}$ independent of $j\in\mathbb Z/t\mathbb Z$, so we have $h = t\cdot \mathrm{hol}(\ell_{(s,t)})$. Therefore, it is enough to compute $\mathrm{hol}(m_t^s)$. We have
\begin{align*}
	 |y_g^{[1]}\cdots y_g^{[t]}x_g^{[1]}\cdots x_g^{[s]}| &\xmapsto{m_t} |y_g^{[1]}\cdots y_g^{[t-1]}(x_g^{[1]}y_g^{[t]} + \omega')x_g^{[2]}\cdots x_g^{[s]}|\\
	&\xmapsto{m_t} \cdots\\
	&\xmapsto{m_t} |y_g^{[1]}\cdots y_g^{[t-1]}x_g^{[1]}\cdots x_g^{[s]}y_g^{[t]} + \sum_{1\leq i\leq s} y_g^{[1]}\cdots y_g^{[t-1]}x_g^{[1]}\cdots x_g^{[i-1]}\omega'x_g^{[i+1]}\cdots x_g^{[s]}|.
\end{align*}
Put $r_{s,t} = \sum_{1\leq i\leq s} x_g^{s-i}y_g^{t-1}x_g^{i-1} \in T(H)$ for $s,t\geq 1$. Since $\irr(|r_{s,t}\omega'|) < +\infty$ and $\omega'$ does not contain any $x_g$'s and $y_g$'s, we can reduce $\irr(|r_{s,t}\omega'|)$ to $0$ by applying $\rho$ to $r_{s,t}$ inside $|r_{s,t}\omega'|$. We put $|r'_{s,t}\omega'| = \rho(|r_{s,t}\omega'|)$; we have $h = t\cdot |r'_{s,t}\omega'|$. In partiular, if $s = t = 1$, we have $r_{s,t} = r'_{s,t} = 1$, which yields $h = |\omega'| = 0$. This gives the trivial relation $0\in|\omid|$, so the case $s = t = 1$ can be discarded.

We write $r'_{s,t}$ in the form $r'_{s,t} = sx_g^{s-1}y_g^{t-1} + b_{s,t}$. Then, $b_{s,t}\in T(H)$ has irregularity $0$ and the number of $x_g$'s and $y_g$'s in each term in $b_{s,t}$ is less than $(s-1)$ and $(t-1)$, respectively.\\[-7pt]
\end{computation}

With this being seen, we define a secondary rewriting rule of monomials (not of a sequence in a word!) for $g\geq 2$ and $s,t\geq 1$, $s+t\geq 3$, by 
\[
	|x_g^{s-1}y_g^{t-1}y_{g-1}x_{g-1}| \xmapsto{\rho_2} |s^{-1}r'_{s,t}(x_{g-1}y_{g-1} + \omega'') - s^{-1}b_{s,t} y_{g-1}x_{g-1}|\,,
\]
where $\omega'' = \sum_{1\leq i<g-1} [x_i,y_i]$. Note that the terms on the right-hand side have irregularity $0$, and $s^{-1}$ is well-defined since we assumed $\mathrm{char}(\mathbb{K})=0$. This rewriting eliminates the monomials $|x_g^{s-1}y_g^{t-1}y_{g-1}x_{g-1}|$ since these do not appear in the right-hand side of $\rho_2$ by the condition $s+t\geq 3$. It always terminates since a polynomial has a finite number of terms.

\begin{remark}
If $s=t=1$, the definition of $\rho_2$ reads $|y_{g-1}x_{g-1}|  \xmapsto{\rho_2} |x_{g-1}y_{g-1}|$, which rewrites nothing at all. Therefore, the monomial $|y_{g-1}x_{g-1}|$ cannot be eliminated.
\end{remark}

Now, we state the main result of this section.

\begin{theorem}\label{thm:basis}
Define the subset $X$, $Y$ of $C$ by 
\begin{align*}
	X &= \{w\in C:\irr(w) = 0,w\neq |x_g^{s-1}y_g^{t-1}y_{g-1}x_{g-1}| \textrm{ for any } s,t \geq 1 \textrm{ with } s+t\geq 3 \},\textrm{ and}\\
	Y &= \{|x_g^sy_g^t|\in C: s,t\geq 1\}. 
\end{align*}
Then, the image of $X\sqcup Y$ in $|T(H)_\omega|$ is a basis for any $g\geq 1$.
\end{theorem}

We need another kind of graph for the proof.
\begin{definition}\ 
\begin{itemize}
	\item A \textit{skew-weighted graph} is a triple $\mathcal{G} = (V(\mathcal{G}), E(\mathcal{G}), \mathsf{wt})$ where $ (V(\mathcal{G}), E(\mathcal{G}))$ is a finite directed graph and $\mathsf{wt}$ is a \textit{weight map} $\mathsf{wt}\colon \mathrm{HE}(\mathcal{G}) \to \mathbb{K}$. Here we denote by $\mathrm{HE}(\mathcal{G}) = \bigsqcup_{w\in V(\mathcal{G})}\mathrm{HE}(w)$ the set of half-edges. We require every vertex of $\mathcal{G}$ to be at least univalent, and $\mathsf{wt}(h) = -\mathsf{wt}(\bar h)$ for each pair of half edges $(h,\bar h)$ comprising an edge. Graphs may have self-loops (i.e. an edge whose half-edges are both located at the same vertex) and multiple edges. For $h\in\mathrm{HE}(\mathcal{G})$, which is a part of a unique edge $\lambda$, we set $\varepsilon_h = +1$ if $h$ is the first half of the directed edge $\lambda$ and $-1$ otherwise.
	\item A \textit{cycle decomposition} of a skew-weighted graph $\mathcal{G}$ is a finite family $(\mathcal{G}_z)_{z\in Z}$ of skew-weighted bivalent subgraphs of $\mathcal{G}$, each with a choice of global orientation (out of possible two of a circle). The weight map $\mathsf{wt}_z$ on $\mathcal{G}_z$ is not necessarily a restriction of $\mathsf{wt}$. Put $\mathrm{HE}(\mathcal{G}_z) = \{h^z_1,\bar h^z_1,\dotsc, h^z_{r_z},\bar h^z_{r_z}\}$ in this cyclic order with respect to the chosen global orientation, where $h^z_i$ and $\bar h^z_i$ is contained in the same edge $\lambda^z_i$. We require that $c_z := \mathsf{wt}_z(h^z_i)$ is independent of $i\in\{1,\dotsc,r_z\}$, $\mathcal{G} = \bigcup_{z\in Z} \mathcal{G}_z$ as a usual graph and 
\begin{align}\label{eq:total}
	\mathsf{wt}(h) = \sum_{z\in Z: \,h\in \mathrm{HE}(\mathcal{G}_z)} \mathsf{wt}_z(h)
\end{align}
for each $h\in\mathrm{HE}(\mathcal{G})$. 
\end{itemize}
\end{definition}

\begin{lemma}\label{lem:cycdecomp}
Let $\mathcal{G}$ be a skew-weighted graph such that $\mathsf{wt}(h)\neq 0$ for all $h\in\mathrm{HE}(\mathcal{G})$ and 
\begin{align}\label{eq:cyclezero}
	\sum_{h\in \mathrm{HE}(w)}\mathsf{wt}(h) = 0
\end{align}
for all $w\in V(\mathcal{G})$. Then, there exists a cyclic decomposition of $\mathcal{G}$.
\end{lemma}
\noindent Proof. If $\mathcal{G}$ has a univalent vertex $w$, there is a unique $h\in\mathrm{HE}(w)$ and we have $\mathsf{wt}(h) = 0$, which contradicts the assumption. Therefore, $\mathcal{G}$ has no univalent vertices.

We proceed by induction on $n = \# E(\mathcal{G})$. If $n=0$, there are no edges. Since we assumed that every vertex is at least univalent, the graph itself is empty. Let $n\geq 1$ and take a bivalent subgraph $\mathcal{G}_1$ in $\mathcal{G}$. If there are no such subgraphs, $\mathcal{G}$ is necessarily a tree, which is impossible since $\mathcal{G}$ has no univalent vertex. Choose a global orientation and put $\mathrm{HE}(\mathcal{G}_1) = \{h_1,\bar h_1,\dotsc, h_{r},\bar h_{r}\}$ in this cyclic order. Put $c = \mathsf{wt}(h_1)$. We define a weight map on $\mathcal{G}_1$ by $\mathsf{wt}_1(h_i) = c$. Then $(\mathcal{G}_1,\mathsf{wt}_1)$ is a skew-weighted graph. We define another weight map on $\mathcal{G}$ by
\[
	\mathsf{wt}_2(h) = \left\{\begin{alignedat}{2} & \mathsf{wt}(h) - \mathsf{wt}_1(h) && \;\textrm{ if } h\in\mathrm{HE}(\mathcal{G}'),\\ &\mathsf{wt}(h) & &\;\textrm{ otherwise.}\end{alignedat}\right.
\]
Let $(\mathcal{G}_2,\mathsf{wt}_2)$ be the skew-weighted subgraph of $\mathcal{G}$ specified by $E(\mathcal{G}_2) = \{\lambda=(h,\bar h)\in E(\mathcal{G}): \mathsf{wt}_2(h)\neq 0\}$. Since $\mathsf{wt}_2(h_1)=0$, we have $\# E(\mathcal{G}_2)<n$. Since $\mathcal{G}_2$ also satisfies the condition \eqref{eq:cyclezero}, we can take a cyclic decomposition $(\mathcal{G}_{2,z})_{z\in Z_2}$ of $\mathcal{G}_2$ by the induction hypothesis. It is readily seen that $(\mathcal{G}_{2,z})_{z\in Z_2}\cup \{\mathcal{G}_1\}$ is a cyclic decomposition. \qed\\

\noindent \textbf{Proof of Theorem \ref{thm:basis}.} The case of $g=1$ is straightforward, so we assume $g\geq 2$. We first show the linear independence of the image of $X$ in $|T(H)_\omega|$, which is equivalent to $\mathrm{Span}_\mathbb{K}(X)\cap |\omid| = 0$ in $|T(H)|$. Take $a\in\mathrm{Span}_\mathbb{K}(X)\cap |\omid|$, which can be expressed as $a = \sum_{\lambda\in\Lambda} c_\lambda (w_\lambda - a_\lambda)$ where $\Lambda$ is a finite set, $c_\lambda\in\mathbb{K}\setminus \{0\}$, $w_\lambda \in C$ with arbitrary irregularity (including $+\infty$) and $a_\lambda\in |T(H)|$ is obtained by applying $\rho$ more than one times to $w_\lambda$. Set $\Lambda_n = \{\lambda\in\Lambda :  \irr(w_\lambda) = n\}$; we have $\irr(w_\lambda) > \irr(a_\lambda)$ for $\lambda \in \Lambda_n$ with $n<+\infty$. For $\lambda\in\Lambda_{+\infty}$, we put $a_\lambda = w'_\lambda + a'_\lambda$ where $w'_\lambda\in C$ with $\irr(w'_\lambda)=+\infty$ and $\irr(a'_\lambda)<+\infty$, and we have an associated path $p_\lambda: w_\lambda \mapsto w'_\lambda$ in $\mathcal{R}$ by remembering how we obtained $a_\lambda$ from $w_\lambda$. Setting $C_{n} = \{w\in C: \irr(w)=n\}$ for $n\in\mathbb{Z}_{\geq 0}\cup \{+\infty\}$, we have a projection
\[
	\mathrm{proj}_{+\infty}\colon |T(H)| \to \mathrm{Span}_\mathbb{K}(C_{+\infty})
\]
induced by the decomposition $C = C_{+\infty} \sqcup (C\setminus C_{+\infty})$. Then, since $\irr(a) = 0$, we have
\begin{align}\label{eq:projinf}
	0 = \mathrm{proj}_{+\infty}(a) &= \mathrm{proj}_{+\infty} \Big(\sum_{\lambda\in\Lambda_{+\infty}} c_\lambda (w_\lambda - a_\lambda)\Big) = \sum_{\lambda\in\Lambda_{+\infty}} c_\lambda (w_\lambda -  w'_\lambda).
\end{align}
Therefore, we have
\begin{align*}
	a =  \sum_{\lambda\in\Lambda\setminus\Lambda_{+\infty}} c_\lambda (w_\lambda - a_\lambda)  +\sum_{\lambda\in\Lambda_{+\infty}} c_\lambda (-a'_\lambda).
\end{align*}
Since $a\in \mathrm{Span}_\mathbb{K}(X)$, we have $\irr(a) = 0$ and
\begin{align*}
	a = \rho(a) &= \rho\Big( \sum_{\lambda\in\Lambda\setminus\Lambda_{+\infty}} c_\lambda (w_\lambda - a_\lambda)  +\sum_{\lambda\in\Lambda_{+\infty}} c_\lambda (-a'_\lambda)\Big)\\
	&= \sum_{\lambda\in\Lambda\setminus\Lambda_{+\infty}} c_\lambda (\rho(w_\lambda) - \rho(a_\lambda))  -\rho\Big(\sum_{\lambda\in\Lambda_{+\infty}} c_\lambda a'_\lambda\Big)\\
	&= -\rho\Big(\sum_{\lambda\in\Lambda_{+\infty}} c_\lambda a'_\lambda\Big).
\end{align*}
Here we used $\rho(w_\lambda) = \rho(a_\lambda)$ for $\lambda\in\Lambda\!\setminus \!\Lambda_{+\infty}$, which follows from Corollary \ref{cor:terminates}.

Now, we construct a skew-weighted graph. Set $V(\mathcal{G}) = \{w_\lambda, w'_\lambda: \lambda\in \Lambda_{+\infty}\}\subset C$ and $E(\mathcal{G}) =  \Lambda_{+\infty}$. Each $\lambda$ corresponds to an edge from $w_\lambda$ to $w'_\lambda$. The weight map is defined by $\mathsf{wt}(h) = \varepsilon_h c_\lambda$ where $h$ is a part of a unique edge $\lambda$ from $w_\lambda$ and $w'_\lambda$. The equality \eqref{eq:projinf} is equivalent to
\begin{align*}
	\sum_{h\in \mathrm{HE}(w)}\mathsf{wt}(h) = 0
\end{align*}
for each $w\in V(\mathcal{G})$. Since $\mathsf{wt}(h)\neq 0$ for all $h$ by the assumption, we can take a cycle decomposition $(\mathcal{G}_z)_{z\in Z}$ of $\mathcal{G}$ by Lemma \ref{lem:cycdecomp}. The right-hand side of \eqref{eq:total} is equal to
\[
	\sum_{z\in Z: \,h=h_i^z} \mathsf{wt}_z(h) + \sum_{z\in Z: \,h=\bar h_i^z} \mathsf{wt}_z(h) = \sum_{z\in Z: \,h=h_i^z} c_z - \sum_{z\in Z: \,h=\bar h_i^z} c_z,
\]
so we obtain
\begin{align*}
	&\sum_{\lambda = (h,\bar h)\in\Lambda_{+\infty}} c_\lambda (w_\lambda -  w'_\lambda - a'_\lambda)\\
	&= \sum_{\lambda = (h,\bar h)\in\Lambda_{+\infty}} \varepsilon_h\mathsf{wt}(h)  (w_\lambda -  w'_\lambda - a'_\lambda)\\
	&= \sum_{\lambda = (h,\bar h)\in\Lambda_{+\infty}} \bigg(\sum_{z\in Z: \,h=h_i^z} \varepsilon_h c_z(w_\lambda -  w'_\lambda - a'_\lambda) - \sum_{z\in Z: \,h=\bar h_i^z} \varepsilon_h c_z (w_\lambda -  w'_\lambda - a'_\lambda)\bigg)\\
	&= \sum_{z\in Z}  \bigg(\sum_{\substack{\lambda = (h,\bar h)\in\Lambda_{+\infty}\\h=h_i^z}} \varepsilon_h c_z(w_\lambda -  w'_\lambda - a'_\lambda) - \sum_{\substack{\lambda = (h,\bar h)\in\Lambda_{+\infty}\\h=\bar h_i^z}} \varepsilon_h c_z (w_\lambda -  w'_\lambda - a'_\lambda)\bigg)\\
	&= \sum_{z\in Z}  c_z \bigg(\sum_{\substack{\lambda = (h,\bar h)\in\Lambda_{+\infty}\\h=h_i^z}} \varepsilon_{h_i^z} (w_{\lambda^z_i} -  w'_{\lambda^z_i} - a'_{\lambda^z_i}) + \sum_{\substack{\lambda = (h,\bar h)\in\Lambda_{+\infty}\\h=\bar h_i^z}} \varepsilon_{h_i^z} (w_{\lambda^z_i} -  w'_{\lambda^z_i} - a'_{\lambda^z_i})\bigg)\\
	&= \sum_{z\in Z}  c_z \sum_{1\leq i\leq r_z} \varepsilon_{h_i^z} (w_{\lambda^z_i} -  w'_{\lambda^z_i} - a'_{\lambda^z_i}) .
\end{align*}
We can choose the global orientation so that $\varepsilon_{h^z_1} = +1$. We set $a'_z = \sum_{1\leq i\leq r_z}  \varepsilon_{h_i^z} a'_{\lambda^z_i} $ for $z\in Z$ so that we have
\[
	 \sum_{\lambda\in\Lambda_{+\infty}} c_\lambda a'_\lambda = \sum_{z\in Z} c_za'_z.
\]
From a bivalent subgraph $\mathcal{G}_z$, we construct a genuine loop $\ell_z$ in $\tilde{\mathcal{R}}^{(s,t)}$ inductively by the following. We start at $w_{\lambda^z_1}$ and traverse along the path $p_{\lambda^z_1}$ to reach $w'_{\lambda^z_1}$. For $i\geq 2$, if the next path $p_{\lambda^z_i}$ is in the correct direction (namely, from $w_{\lambda^z_i}$ to $w'_{\lambda^z_i}$, in which case $\varepsilon_{h_i^z} = +1$), traverse $p_{\lambda^z_i}$; if not, traverse along an arbitrary directed path $p'_{\lambda^z_i}$ from $w_{\lambda^z_i}$ to $w'_{\lambda^z_i}$, whose existence is guaranteed by Lemma \ref{lem:dirconn}. Then, the holonomy along $\ell_z$ is equal to a scalar multiple of $|r'_{s,t}\omega'|$ by Corollary \ref{cor:deform} and Computation \ref{comp:hol}. For $\lambda^z_i$ such that $\varepsilon_{h^z_i} = -1$, we have a loop $p_{\lambda^z_i}p'_{\lambda^z_i}$ in $\mathcal{R}$, whose holonomy is given by $\mathrm{hol}(p_{\lambda^z_i}) + \mathrm{hol}(p'_{\lambda^z_i}) = q_{z,i}|r'_{s_z,t_z}\omega'|$ for some $q_{z,i}\in\mathbb{K}$. Then, we have
\begin{align*}
	a'_z - \mathrm{hol}(\ell_z) &= \sum_{1\leq i\leq r_z}  \varepsilon_{h^z_i} \mathrm{hol}(p_{\lambda^z_i}) - \sum_{\substack{1\leq i\leq r_z\\\varepsilon_{h^z_i} = +1}}  \mathrm{hol}(p_{\lambda^z_i}) - \sum_{\substack{1\leq i\leq r_z\\\varepsilon_{h^z_i} = -1}}  \mathrm{hol}(p'_{\lambda^z_i}) \\
	&= \sum_{\substack{1\leq i\leq r_z\\\varepsilon_{h^z_i} = -1}} \varepsilon_{h^z_i} \mathrm{hol}(p_{\lambda^z_i}) -  \mathrm{hol}(p'_{\lambda^z_i})\\
	&= \sum_{\substack{1\leq i\leq r_z\\\varepsilon_{\lambda_i} = -1}} -q_{z,i} |r'_{s_z,t_z}\omega'|\,.
\end{align*}
Therefore, we have $\rho(a'_z) = q_z|r'_{s_z,t_z}\omega'|$ for some $q_z\in\mathbb{K}$. Hence, we obtain
\[
	a =  -\rho\Big(\sum_{\lambda\in\Lambda_{+\infty}} c_\lambda a'_\lambda\Big) = -\sum_{z\in Z} c_z \rho(a'_z) = -\sum_{z\in Z} c_z q_z|r'_{s_z,t_z}\omega'| = - \sum_{s,t\geq 1} c_{s,t} |r'_{s,t}\omega'| = -  \sum_{\substack{s,t\geq 1\\s+t\geq 3}}c_{s,t} |r'_{s,t}\omega'|
\]
for some $c_{s,t} \in\mathbb{K}$. We introduce another projection
\[
	\mathrm{proj}'\colon \mathrm{Span}_\mathbb{K}(C_0) \to \mathrm{Span}_\mathbb{K}(C_0\setminus X)
\]
induced from the decomposition $C_0 = X \sqcup (C_0\setminus X)$. Since $a\in \mathrm{Span}_\mathbb{K}(X)$, We have
\[
	0 = \mathrm{proj}'(a) = \sum_{\substack{s,t\geq 1\\s+t\geq 3}} c_{s,t}\cdot s |x_g^{s-1}y_g^{t-1}y_{g-1}x_{g-1}|.
\]
The monomials $|x_g^{s-1}y_g^{t-1}y_{g-1}x_{g-1}|$ are distinct, so we obtain $c_{s,t}\cdot s = 0$ for all $s,t\geq 1$ with $s+t\geq 3$, which is equivalent to $c_{s,t} = 0$. Therefore, we obtain $a=0\,$; this shows the linear independency of the image of $X$ in $|T(H)_\omega|$.

Next, considering the $\mathbb{K}$-algebra map
\[
	\phi\colon T(H)_\omega \to \mathbb{K}[x_g,y_g]\colon x_i,y_i \mapsto 0\, (i<g), \;x_g\mapsto x_g,\; y_g\mapsto y_g\,,
\]
the composite
\[
	\mathrm{Span}_\mathbb{K}(X\sqcup Y) \to |T(H)_\omega| \xrightarrow{|\phi|}\mathbb{K}[x_g,y_g]\twoheadrightarrow \mathbb{K}[x_g,y_g]/\langle x_g^s, y_g^t\colon s,t\geq 0\rangle_{\mathbb{K}\textrm{-vect}}
\]
sends $X$ to 0 and maps $\mathrm{Span}_\mathbb{K}(Y)$ to the image isomorphically. Therefore, the image of $X\sqcup Y$ in $|T(H)_\omega|$ is linearly independent.

Finally, we show $\mathrm{Span}_\mathbb{K}(X\sqcup Y)$ surjects onto $|T(H)_\omega|$. For $a\in |T(H)|$, apply $\rho$ to the terms in $a$ with $\irr = +\infty$ so that they are contained in $\mathrm{Span}_\mathbb{K}(Y)$. Denote the (not necessarily unique) result by $a'$. Next, apply $\rho$ to the terms in $a'$ with $\irr < +\infty$ until their irregularity becomes zero. Denoting the result by $a''$, apply $\rho_2$ to all the terms in $a''$ which are scalar multiples of $|x_g^{s-1}y_g^{t-1}y_{g-1}x_{g-1}|$ for $s+t\geq 3$ to eliminate them. The final result $a'''$ is now contained in $\mathrm{Span}_\mathbb{K}(X\sqcup Y)$, while $a', a', a''$ and $a'''$ are all equal modulo $|\omid|$. \qed\\

\section{The Closed Surface Case}\label{sec:closed}
In this section, we solve the formality problem for a closed surface $\Sigma = \Sigma_{g,0}$. The following is already known in \cite{akkn}.

\begin{theorem}
The completed Goldman--Turaev Lie bialgebra $|\widehat{\mathbb{K}\pi}/\mathbb{K}1|$ for a closed surface is formal.
\end{theorem}
\noindent Proof. This is a consequence of the formality on a surface with only one boundary component. First of all, since a framing does not matter in this case, the set of formality isomorphisms is non-empty by Theorem 7.1 of \cite{akkn} for any $g\geq 1$. Therefore, we can take a formality isomorphism $F\colon \widehat{\mathbb{K}\pi_1(\Sigma_{g,1})} \to \mathrm{gr}(\widehat{\mathbb{K}\pi_1(\Sigma_{g,1})})$. This descends to $F\colon \widehat{\mathbb{K}\pi_1(\Sigma_{g,0})} \to \mathrm{gr}(\widehat{\mathbb{K}\pi_1(\Sigma_{g,0})})$ as being the formality isomorphism forces $F$ to satisfy the condition $F(\xi) = \omega$. Since the Lie bialgebra structure of $|\widehat{\mathbb{K}\pi_1(\Sigma_{g,0})}/\mathbb{K}1|$ is the natural quotient of $|\widehat{\mathbb{K}\pi_1(\Sigma_{g,1})}|$, the claim follows.\qed\\

The theorem above states that the set of formality isomorphisms is a torsor over their pro-unipotent automorphism groups. At this point, describing the automorphism group suffices. 

We are done with linear categories and return to the usual algebras. Let $\hat L(H)_\omega = \hat L(H)/\langle\omega\rangle$ be the quotient Lie algebra by the ideal generated by $\omega$, and $\hat T(H)_\omega = U\hat L(H)_\omega$, which is a quotient of $\hat T(H)$ by the (complete) two-sided ideal generated by $\omega$. Now, we take the standard resolution of $\mathbb{K}$ as a left $\hat T(H)_\omega$-module
\begin{align*}
	0 \to P_1 \xrightarrow{\partial_1} &P_0 \xrightarrow{\partial_0} \hat T(H)_\omega \xrightarrow{\mathrm{aug}} \mathbb{K}\to 0
\end{align*}
where $P_0 = \hat T(H)_\omega\otimes H$, $P_1 = \hat T(H)_\omega\otimes \mathbb{K}\omega$, and
\[
	\partial_1(1\otimes \omega) = \sum_i x_i\otimes y_i - y_i\otimes x_i\,,\;\partial_0(a\otimes v) = av\,.
\]
The proof of exactness can be obtained by taking the associated graded of the proof of Lemma 8.4 in \cite{toyo}. Now we define connections on each piece:
\begin{align*}
	\nabla'_{0,H}(1\otimes x_i) = \nabla'_{0,H}(1\otimes y_i) = 0 \textrm{ on } P_0\textrm{ and } \nabla'_{1,H}(1\otimes \omega) = 0 \textrm{ on } P_1\,.
\end{align*}
Put $\nabla'_{\bullet,H} = \{\nabla'_{i,H}\}_{i\geq 0}$; this comprises a flat \textit{homological connection} defined in \cite{toyo}, Section 7.\\[-7pt]

The following is the main result of this paper.

\begin{theorem}\label{thm:kvclosed}
The pro-unipotent part of the automorphism group of the associated graded of the Goldman--Turaev Lie bialgebra $(|\hat T(H)_\omega/\mathbb{K}1|,[\cdot,\cdot]_\mathrm{gr},\delta_\mathrm{gr})$ is given by $\mathrm{KRV}_{(g,0)} = \Exp(\mathfrak{krv}_{(g,0)})$, where
\[
	\mathfrak{krv}_{(g,0)} := \{g\in\Der^+(\hat L(H)_\omega)\colon \sdiv^{\nabla'_{\bullet,H}}(g) \in \Ker(|\bar\Delta_\omega|)\}\,.
\]
Here, $\bar\Delta_\omega\colon \hat T(H)_\omega \to \hat T(H)_\omega\hat \otimes \hat T(H)_\omega$ is the reduced coproduct: $\bar\Delta_\omega(x) = \Delta_\omega(x) - x\otimes 1 - 1\otimes x$.
\end{theorem}

We need some lemmas for the proof. The symbol $[\cdot,\cdot]$ denotes the usual commutator, while $[\cdot,\cdot]\gr$ denotes the associated graded of the Goldman bracket.

\begin{lemma}\label{lem:comm}
Fix a symplectic basis $(x_i,y_i)_{1\leq i\leq g}$ of $H$ as before. Let $z\in T(H)_\omega$ be a primitive non-zero monomial in $(x_i,y_i)_{1\leq i\leq g}$ with $\deg(z)\leq 2$ (i.e., not of the form $(z')^m$ for $m\geq 2$), and $b\in T(H)_\omega$ a homogeneous element. If $|bz^\ell|=0$ for a sufficiently large $\ell$, we have $b\in [z, T(H)_\omega]$ and therefore $|b| = 0$.
\end{lemma}
\noindent Proof. This is a version of Lemma A.3 in \cite{akkn}, so we basically follow them. If $g=1$, the spaces $T(H)_\omega$ and $|T(H)_\omega|$ are both isomorphic to the (commutative) polynomial algebra $S(H)$ so we have $b=0$. Now assume $g\geq 2$. We prove by induction on $\deg(b)$. If $\deg(b)=0$, we have $b=0$ since $b$ is just a scalar. Suppose $\deg(b)\geq 1$.

\underline{Case of $\deg(z)=1$}. We can assume $z = x_1$ by applying a graded algebra automorphism of $T(H)_\omega$ that sends $z$ to $x_1$. Let $\tilde b\in T(H)$ be a lift of $b$ such that there is no sequence $y_gx_g$ appearing in $\tilde b$. Then, we have $\irr(|\tilde bx_1^\ell|)=0$ by the assumption $g\geq 2$. Furthermore, each term in $|\tilde bx_1^\ell|$ is not of the form $|x_g^{s-1}y_g^{t-1}y_{g-1}x_{g-1}|$ for any $s,t\geq 1$ since $\ell$ is sufficiently large. Then, we have $|\tilde bx_1^\ell|\in\mathrm{Span}_\mathbb{K}(X)$. By Theorem \ref{thm:basis}, we obtain $|\tilde bx_1^\ell| = 0$ in $|T(H)|$.

\underline{Case of $\deg(z)=2$}. We can assume $z = x_1x_2$ or $x_1y_1$ for the same reason as above. Let $\tilde b\in T(H)$ be a lift of $b$ such that there is no sequence $y_gx_g$ appearing in $\tilde b$. Then, we have $\irr(|\tilde b z^\ell|) = 0$ by the assumption $g\geq 2$. Furthermore, each term in $|\tilde bz^\ell|$ is not of the form $|x_g^{s-1}y_g^{t-1}y_{g-1}x_{g-1}|$ for any $s,t\geq 1$ since $\ell$ is sufficiently large. Then, we have $|\tilde bz^\ell|\in\mathrm{Span}_\mathbb{K}(X)$. By Theorem \ref{thm:basis}, we obtain $|\tilde bz^\ell| = 0$ in $|T(H)|$.

Therefore, in any case, we have $|\tilde bz^\ell| = 0$ in $|T(H)|$. Now uniquely express $\tilde b$ as $\tilde b = z\tilde b_1 + \tilde b_2z + \tilde b_3$ where $\tilde b_i\in T(H)$, the terms in $\tilde b_2$ cannot divided by $z$ from left, and the terms in $\tilde b_3$ cannot be divided by $z$ from neither side. We have
\[
	0 = |\tilde b_1z^{\ell+1} + \tilde b_2z^{\ell+1} + \tilde b_3z^\ell|.
\]
The first and second terms have $(\ell + 1)$ consecutive $z$'s in their terms, while none of $|\tilde b_3z^\ell|$ does by the assumption. Therefore we have $|\tilde b_3z^\ell| = 0$ and hence $|(\tilde b_1+\tilde b_2)z^{\ell+1}|=0$. Since $\ell$ is sufficiently large and $\tilde b_3$ cannot be divided by $z$ from neither side, we have $\tilde b_3=0$. In addition, by the induction hypothesis, we have $\tilde b_1 + \tilde b_2 = [z, \tilde b_4]$ for some $\tilde b_4\in T(H)$. We obtain
\[
	\tilde b = z\tilde b_1 + ([z,\tilde b_4] - \tilde b_1)z = [z,\tilde b_1 + \tilde b_4 z].
\]
Therefore, we have $b = [z,b_1 + b_4z]$; this completes the proof.\qed\\[-7pt]

\begin{lemma}\label{lem:wedgetwo}
Putting $\mathfrak{g} = |\hat T(H)_\omega/\mathbb{K}1|$ regarded as a Lie algebra, we have $(\mathfrak{g}\otimes \mathfrak{g})^\mathfrak{g} = 0$.
\end{lemma}
\noindent Proof. Since the defining equations of the space $(\mathfrak{g}\otimes \mathfrak{g})^\mathfrak{g}$ are linear, we can assume $\mathbb{K}$ to be algebraically closed. Furthermore, since the space $(\mathfrak{g}\otimes \mathfrak{g})^\mathfrak{g}$ is graded, we can work in the non-completed $T(H)_\omega$ without any problem. Then, it is reduced to Crawley--Boevey's result that the quiver variety is generically symplectic: let $Q$ be the quiver with one vertex and $g$ arrows. Then, the \textit{preprojective algebra} $\Pi$ associated with $Q$ is isomorphic to $T(H)_\omega$ as a $\mathbb{K}$-algebra. Since we know the Lie algebra structure on $|T(H)_\omega|$ coincides with the necklace Lie bracket on $|\Pi|$ (see Section 6 of \cite{gtjohnson}, for example), Theorem 8.6.1 (ii) of \cite{quiver} implies $Z(|T(H)_\omega|) = Z(|\Pi|) = |\mathbb{K}1|$.

Now take a homogeneous element $x$ of $(\mathfrak g\otimes \mathfrak g)^\mathfrak{g}$ with $\deg(x)=r$. Denote by $\tilde x\in |T(H)_\omega|^{(\geq 1)\,\otimes 2}$ the unique lift of $x$ along the natural projection $|\hat T(H)_\omega|^{\otimes 2}\twoheadrightarrow\mathfrak g^{\otimes 2}$. The condition $x\in (\mathfrak g\otimes \mathfrak g)^\mathfrak{g}$ is equivalent to $[\tilde y, \tilde x]\gr \in |1|\otimes |T(H)_\omega| +  |T(H)_\omega| \otimes |1|$ for all homogeneous $y\in \mathfrak{g}$ and its unique lift $\tilde y\in|T(H)_\omega|^{(\geq 1)}$. Put $\tilde x = \sum_{0<i<r} \sum_{\lambda\in \Lambda_i} \tilde x_{i,\lambda}' \otimes w_{i,\lambda}$ where $\tilde x_{i,\lambda}' \in |T(H)_\omega|$ with $\deg(x_{i,\lambda}') = i$ and $w_{i,\lambda}$ is an element of some fixed homogeneous basis of $|T(H)_\omega|$. We have
\[
	[\tilde y, \tilde x]\gr  = \sum_{0<i<r} \sum_{\lambda\in \Lambda_i} ([\tilde y,\tilde x_{i,\lambda}']\gr \otimes w_{i,\lambda} + \tilde x_{i,\lambda}' \otimes [\tilde y,w_{i,\lambda} ]\gr ).
\]
If $\deg(y)\geq 2$, we have $\deg([\tilde y, \tilde x_{i,\lambda}']\gr) = \deg(y) + i - 2 \geq i \geq 1$ and, similarly, $\deg([\tilde y, w_{i.\lambda}]\gr) \geq 1$ by Proposition \ref{prop:filtdeg}. Therefore, in this case, we obtain $[\tilde y, \tilde x]\gr = 0$ for any $y$ with $\deg(y)\geq 2$.

Suppose $\tilde x \neq 0$. Take the largest $i$ such that $\tilde x_{i,\lambda}'\neq 0$ for some $\lambda$. Then, if $\deg(y)\geq 3$, we have $[\tilde y, \tilde x_{i,\lambda}']\gr = 0$ for all $\lambda\in\Lambda_i$ since the map $[\tilde y,\cdot\,]\gr$ have degree $\geq 1$ and we took $i$ to be the largest. Let $u_{i,\lambda}'\in\Der_\mathbb{K}(T(H)_\omega)$ be a representative of $\sigma\gr(\tilde x_{i,\lambda}')\in\mathrm{HH}^1(\Pi)$. Then, $[\tilde y, \tilde x_{i,\lambda}']\gr = 0$ is equivalent to $|u_{i,\lambda}'(\tilde y)| = 0$. In particular, for any monomial $a\in\mathfrak g$ and $\ell\geq 3$ we have
\begin{align*}
	0 = |u_{i,\lambda}'(\tilde a^\ell)| = \ell |u_{i,\lambda}'(\tilde a) \tilde a^{\ell-1}|.
\end{align*}
By Lemma \ref{lem:comm}, for $a\in\mathfrak g$ with $\deg(a)\leq 2$, we have $|u_{i,\lambda}'(\tilde a)| = 0$. In fact, if $\tilde a$ is primitive, we can directly apply Lemma \ref{lem:comm}. If $\tilde a = |z^2|$ for some $z\in H$, we have $u'_{i,\lambda}(z) = [z, b]$ for some $b\in T(H)_\omega$ again by Lemma \ref{lem:comm} and hence
\[	
	|u'_{i,\lambda}(z^2)| = 2|zu'_{i,\lambda}(z)| = 2|z[z,b]| = 0. 
\]
In conclusion, we have $|u_{i,\lambda}'(\tilde y)| = 0$ for any $y\in\mathfrak g$, which, in turn, implies $\tilde x_{i,\lambda}'$ is in the centre $Z(|T(H)_\omega|)$. This is a contradiction since $\deg(x'_{i,\lambda})\geq 1$ and $\tilde x'_{i,\lambda}\neq 0$ for some $\lambda$. Therefore, we have $\tilde x = 0$ and hence $x=0$. \qed\\[-7pt]

\begin{lemma}\label{lem:kviclosed}
Any $F\in\Aut^+_\mathrm{Hopf}(\hat T(H)_\omega)$ preserves the graded Hamiltonian flow $\sigma_\mathrm{gr}\colon |\hat T(H)_\omega| \to \mathrm{HH}^1(\hat T(H)_\omega)$.
\end{lemma}
\noindent Proof. We have $F = \exp(f)$ for some $f\in\Der^+(\hat L(H)_\omega)$. We will show that any element $\Der^+(\hat L(H)_\omega)$ can be lifted to $\Der_\omega^+(\hat L(H)) = \{\tilde f\in\Der^+(\hat L(H)):\tilde f(\omega)=0\}$. Assume $f$ is homogeneous and take any lift $\tilde f$ of $f$ to $\Der^+(\hat L(H))$. We are done if $\tilde f(\omega)=0$, but it is not the case in general; instead, $\tilde f(\omega)$ is an element of the ideal $[H,[H,\cdots [H,\omega]\cdots]]$ generated by $\omega$. Since $\deg(f(\omega))\geq 1+ \deg(\omega)=3$, we can write $\tilde f(\omega) = \sum_i [x_i,a_i] + [y_i,b_i]$ with $a_i, b_i\in [H,[H,\cdots [H,\omega]\cdots]]$. Now define $f' \in \Der^+(\hat L(H))$ by $f'(x_i) = \tilde f(x_i) + b_i$ and $f'(y_i) = \tilde f(y_i) - a_i$. Then, $f'$ is also a lift of $f$, and 
\begin{align*}
	f'(\omega) &= \sum_i [f'(x_i),y_i] + [x_i,f'(y_i)]\\
	&= \tilde f(\omega) + \sum_i [ b_i,y_i] + [x_i, -a_i]\\
	&= \sum_i [x_i,a_i] + [y_i,b_i] + [b_i,y_i] + [x_i, -a_i]\\
	&= 0.
\end{align*}
Therefore, $F\in\Aut^+_\mathrm{Hopf}(\hat T(H)_\omega)$ can also be lifted to an element of $\Aut^+_{\mathrm{Hopf},\omega}(\hat T(H))$. Now that Remark 6.10 of \cite{akkn} says $\Aut^+_{\mathrm{Hopf},\omega}(\hat T(H))$ preserves $\sigma\gr$, the claim immediately follows.\qed\\[-7pt]

\begin{remark}
The above lemma also follows from the fact that $\sigma_\mathrm{gr}$ is a composition
\[
	|\hat T(H)_\omega|= \mathrm{HH}_0(\hat T(H)_\omega) \xrightarrow{\mathsf{B}} \mathrm{HH}_1(\hat T(H)_\omega) \cong \mathrm{HH}^1(\hat T(H)_\omega)
\]
where $\mathsf{B}$ is Connes' differential, and the isomorphism is given by the Poincar\'e duality. Since $\mathsf{B}$ is functorial, it commutes with $F$. The Poincar\'e duality is given by the cap product with the fundamental class in $\mathrm{HH}_2(\hat T(H)_\omega)$, which is preserved by $F$ since we have assumed  $\mathrm{gr}(F) = \id$. \\[-7pt]
\end{remark}

The connection $\nabla'\!_{\bullet,H}$ induces a flat homological connection $\nabla\!_{\bullet,H}$ on $0 \to Q_1\to Q_0 \to  \Omega^1\hat T(H)_\omega \to 0$, where the resolution is given by 
\[
	Q_0 = \hat T(H)_\omega\otimes H \otimes \hat T(H)_\omega \;\textrm{ and }\; Q_1 =  \hat T(H)_\omega\otimes \mathbb{K}\omega \otimes \hat T(H)_\omega.
\]	

\begin{lemma}\label{lem:cobgrconn}
The associated graded of the Turaev cobracket is equal to the composition $\Div^{\nabla_{\bullet,H}}\circ \,\sigma_\mathrm{gr}$.
\end{lemma}
\noindent Proof. This is obtained by taking the associated graded of Theorem 8.2 of \cite{toyo}. \qed\\

\noindent\textbf{Proof of Theorem \ref{thm:kvclosed}.} We basically imitate the proof in \cite{akkn}, which is for the surface with non-empty boundary.

Suppose that $G\in\Aut^+(\hat L(H)_\omega)$ induces an automorphism of the graded Goldman--Turaev Lie bialgebra. Since $G$ automatically preserves the Lie bracket by Lemma \ref{lem:kviclosed}, this is equivalent to that $G$ only preserves the Lie cobracket. By Lemma \ref{lem:cobgrconn}, this is equivalent to
\[
	G\circ \Div^{\nabla_{\bullet,H}}\circ\, \sigma\gr-  \Div^{\nabla_{\bullet,H}}\circ\,\sigma\gr \circ G = 0 .
\]
Since $G\in\Aut^+(\hat L(H)_\omega)$, we can put $g = \log(G) \in \Der^+(\hat L(H)_\omega)$. Then, the above is further equivalent to 
\[
	g\circ \Div^{\nabla_{\bullet,H}}\circ\, \sigma\gr  -  \Div^{\nabla_{\bullet,H}}\circ\, \sigma\gr\circ g= 0 .
\]
By Lemma \ref{lem:kviclosed}, we have $\Ad_G\circ\,\sigma_\mathrm{gr}=  \sigma_\mathrm{gr}\circ G$ on $\mathfrak g$, which in turn implies $\ad_g\circ \,\sigma_\mathrm{gr} =  \sigma_\mathrm{gr}\circ g$. Therefore, we have
\[
	g\circ \Div^{\nabla_{\bullet,H}}\circ\, \sigma\gr  -  \Div^{\nabla_{\bullet,H}}\circ\ad_g \circ\, \sigma\gr= 0 .
\]
By the definition of the infinitesimal adjoint action, we have $g\circ \Div^{\nabla_{\bullet,H}}  -  \Div^{\nabla_{\bullet,H}}\circ\ad_g = \Tr(\ad_g\!\nabla_{\bullet,H})$, which is equal to $d(\Div^{\nabla_{\bullet,H}}(g))$ by Corollary A.5 from \cite{toyo3}.
Hence, for any $a\in\mathfrak g$ and $u = \sigma\gr(a)$, we have
\[
	0 = u(\Div^{\nabla_{\bullet,H}}(g)) = [a,\Div^{\nabla_{\bullet,H}}(g)]\gr.
\]
Therefore, $G$ being a formality isomorphism is equivalent to $\Div^{\nabla_{\bullet,H}}(g) \in (\mathfrak{g}\otimes \mathfrak{g})^\mathfrak{g}$. By Lemma \ref{lem:wedgetwo}, this is equivalent to $\Div^{\nabla_{\bullet,H}}(g) = 0$ in $\mathfrak{g}^{\otimes 2}$.

Now we have $\Div^{\nabla_{\bullet,H}}(g) = |\tilde\Delta_\omega|(\sdiv^{\nabla'_{\bullet,H}}(g))$ by Section 6 of \cite{toyo}. Since the antipode is bijective, $\Div^{\nabla_{\bullet,H}}(g) = 0$ is equivalent to $|\Delta_\omega|(\sdiv^{\nabla'_{\bullet,H}}(g)) = 0\in\mathfrak g^{\otimes 2}$. In addition, $\sdiv^{\nabla'_{\bullet, H}}(g)$ has no constant term since the connection and hence the divergence both have degree 0 while $g$ has a positive degree. Therefore, using a canonical map $\mathfrak g \xrightarrow{\cong} |\hat T(H)_\omega|^{(\geq 1)} \subset |\hat T(H)_\omega|$, the equality is restated as 
\[
	 |\Delta_\omega|(\sdiv^{\nabla'_{\bullet,H}}(g)) - |1|\otimes \sdiv^{\nabla'_{\bullet,H}}(g) - \sdiv^{\nabla'_{\bullet,H}}(g)\otimes |1| = 0\;\textrm{ in } |\hat T(H)_\omega|^{\otimes 2}.
\]
This is exactly the statement of the theorem.\qed\\

\section{The Kernel of the Reduced Coproduct}\label{sec:redcoprod}

In this section, we determine the space $\Ker(|\bar\Delta_\omega|)$ appeared in Theorem \ref{thm:kvclosed} in low degrees. Let $L(H)$ be the free Lie algebra over $H$ (without completion) and $L(H)_\omega$ the quotient of $L(H)$ by the ideal generated by $\omega$. We have $T(H)_\omega = UL(H)_\omega$.\\

\noindent\textbf{Notations.} The superscript ${}^{(d)}$ denotes the degree $d$ part for $d\in \mathbb Z$. For $\mathbb{K}$-subspaces $Z_1, Z_2$ of a $\mathbb{K}$-algebra $A$, we denote by $|Z_1Z_2|$ the image of the composite
\[
	Z_1\otimes Z_2 \xhookrightarrow{\mathrm{incl}^{\otimes 2}} A\otimes A \xrightarrow{\mathrm{mult}} A\twoheadrightarrow |A|.
\]

The main result of this section is the following.

\begin{proposition}\label{prop:ckom}
In low degrees, $\Ker(|\bar\Delta_\omega|\colon|T(H)_\omega|\to |T(H)_\omega|^{\otimes 2})$ is given explicitly by
\begin{align*}
	\Ker(|\bar\Delta_\omega|)^{(1)} &= H,\\
	\Ker(|\bar\Delta_\omega|)^{(2)} &= 0,\\
	\Ker(|\bar\Delta_\omega|)^{(3)} &=|HL(H)_\omega^{(2)}|\cong  \wedge^3 H/|H\omega|,
\end{align*}
for any $g\geq 1$, and $\Ker(|\bar\Delta_\omega|)^{(4)} = |HL(H)_\omega^{(3)}|$ for $g \neq 2$.
\end{proposition}
\begin{remark}
For $g=2$, we have a strict inclusion $\Ker(|\bar\Delta_\omega|)^{(4)} \supset |HL(H)_\omega^{(3)}|$.
\end{remark}

\subsection{The Case of a Free Associative Algebra}

We will start with the case of the free associative algebra $T(H)$ and later apply the result to the case of the quotient algebra $T(H)_\omega$.

\begin{lemma}\label{lem:ckth}
In low degrees, $\Ker(|\bar\Delta|\colon|T(H)|\to |T(H)|^{\otimes 2})$ is given explicitly by
\begin{align*}
	\Ker(|\bar\Delta|)^{(1)} &= H,\\
	\Ker(|\bar\Delta|)^{(2)} &= 0,\\
	\Ker(|\bar\Delta|)^{(3)} &= |HL(H)^{(2)}|\cong \wedge^3 H,\textrm{ and}\\
	\Ker(|\bar\Delta|)^{(4)} &= |HL(H)^{(3)}|
\end{align*}
for any $g\geq 1$.
\end{lemma}
\noindent Proof. The degree 1 part is easily checked. In degree 2, the reduced coproduct reads
\[
	\bar\Delta\colon S^2(H)=|HH| \to H\otimes H\colon |xy| \mapsto |x|\otimes |y| + |y|\otimes |x|.
\]
This is an isomorphism onto the image.

Next, the inclusion $\Ker(|\bar\Delta|)^{(d)} \supset |HL(H)^{(d-1)}|$ for $d\geq 3$ is checked by the following calculation: for $x\in H$ and $\ell\in L$,
\[
	|\Delta(x\ell)| = |\Delta(x)\Delta(\ell)| = |(x\otimes 1 + 1\otimes x)(\ell\otimes 1 + 1\otimes \ell)| = |x\ell\otimes 1 + x\otimes \ell + \ell\otimes x + 1\otimes x\ell|.
\]
Since $|L(H)^{(d-1)}| = 0$ for $d\geq 3$, $|x\otimes \ell + \ell\otimes x| = |x|\otimes |\ell| + |\ell|\otimes |x| $ vanishes if $\deg(\ell)\geq 2$.

Considering the multi-grading, we show the converse is true for $d = 3, 4$. Take a $\mathbb{K}$-basis $(v_i)_{1\leq i\leq 2g}$ of $H$. The multi-degree is given by the decomposition $T(H) = \bigoplus_{\mathbf{d}\in\mathbb{Z}^{2g}} T(H)_{\mathbf d}$ specified by $\deg(v_i) = e_i$. Here $e_i$ is the $i$-th standard unit vector of $\mathbb{Z}^{2g}$. The homogeneous part $T(H)_{\mathbf d}$ is zero if any coefficient of $e_i$ is negative. The coproduct $\Delta$ and hence the reduced $\bar\Delta$ are degree $0$ with respect to this grading, so we only have to compute the kernel for each fixed multi-degree. The total degree map is defined by $|\cdot|\colon \mathbb{Z}^{2g}\to \mathbb{Z}\colon e_i \mapsto 1$.

In degree 3, only possible multi-degrees are, up to permutations of $e_i$'s,
\begin{enumerate}[(i)]
	\item $(3,0,\dotsc,0)$;
	\item $(2,1,0,\dotsc,0)$; or
	\item $(1,1,1,0,\dotsc,0)$.
\end{enumerate}
We compute the kernel in each case and show that it is contained in $|HL(H)^{(2)}|$.
\begin{enumerate}[(i)]
	\item The only mononial with multi-degree $(3,0,\dotsc,0)$ is $|v_1^3|$. Then, we have $0=\bar\Delta(c|v_1^3|) = 3c|v_1|\otimes |v_1^2|$ for $c\in\mathbb{K}$, which implies $c=0$.
	\item Similarly, for $c\in\mathbb{K}$, we have $0=\bar\Delta(c|v_1^2v_2|) = 2c|v_1|\otimes |v_1v_2| + c|v_2|\otimes |v_1^2|$. Again $c=0$. 
	\item For $c,c'\in\mathbb{K}$, we have
	\begin{align*}
		0 &=\bar\Delta(c|v_1v_2v_3| + c'|v_1v_3v_2|)\\
		&= c(|v_1|\otimes |v_2v_3| + |v_2|\otimes |v_1v_3| + |v_3|\otimes |v_1v_2|) + c'(|v_1|\otimes |v_3v_2| + |v_3|\otimes |v_1v_2| + |v_2|\otimes |v_1v_3|).
	\end{align*}
	Therefore, we obtain $c+c'=0$, and we have $c|v_1v_2v_3| + c'|v_1v_3v_2| = c|v_1[v_2,v_3]|\in|HL(H)^{(2)}|$.
\end{enumerate}
This completes the case of degree 3. \\
	
Now, we deal with the degree $4$ case. Only possible multi-degrees with $|\mathbf d| = 4$ are, up to permutations of $e_i$'s, one of the following:
\begin{enumerate}[(i)]
	\item $(4,0,\dotsc,0)$;
	\item $(3,1,0,\dotsc,0)$;
	\item $(2,2,0,\dotsc,0)$;
	\item $(2,1,1,0,\dotsc,0)$; or
	\item $(1,1,1,1,0,\dotsc,0)$.
\end{enumerate}
We compute the kernel in each case and show that it is contained in $|HL(H)^{(3)}|$.
\begin{enumerate}[(i)]
	\item The only mononial with multi-degree $(4,0,\dotsc,0)$ is $v_1^4$. We put $a = cv_1^4$ with $c\in\mathbb{K}$. Then, we have
	\[
		0 = |4c\,v_1\otimes v_1^3 + 6c\,v_1^2\otimes v_1^2|.
	\]
	by seeing the degree $(1,3)$ part of $|\bar\Delta(a)|$. Therefore, $|a|\in \Ker(|\bar\Delta|)$ implies $c=0$.
	\item In multi-degree $(3,1,0,\dotsc,0)$, it is enough to consider $a = cv_1^3v_2$. Then, 
	\[
		0 = |3cv_1\otimes v_1^2v_2 + cv_2\otimes v_1^3|.
	\]
	Therefore, $|a|\in \Ker(|\bar\Delta|)$ implies $c=0$.
	\item In multi-degree $(2,2,0,\dotsc,0)$, it is enough to consider the case $a = c_1v_1^2v_2^2 + c_2v_1v_2v_1v_2$. Then, 
	\[
		0 = |2c_1(v_1\otimes v_1v_2^2 + v_2\otimes v_1^2v_2) + 2c_2(v_1\otimes v_1v_2^2 + v_2\otimes v_1^2v_2)|.
	\]
	Therefore, $|a|\in \Ker(|\bar\Delta|)$ implies $c_1 + c_2=0$. We have
	\begin{align*}
		|a| &= c_1|v_1(v_1v_2^2  - v_2v_1v_2)|\\
		&= c_1|v_1[v_1,v_2]v_2|\\
		&= \frac{c_1}{2}|v_1[[v_1,v_2],v_2]|.
	\end{align*}
	Thus we have $|a|\in |HL(H)^{(3)}|$.
	\item In multi-degree $(2,1,1,0,\dotsc,0)$, it is enough to consider the case $a = c_1v_1^2v_2v_3 + c_2v_1v_2v_1v_3 + c_3v_2v_1^2v_3$. Then, 
	\begin{align*}
		0 &= |c_1(2v_1\otimes v_1v_2v_3 + v_2\otimes v_1^2v_3 + v_3\otimes v_1^2v_2)\\
		&\qquad + c_2(v_1\otimes v_2v_1v_3 + v_2\otimes v_1^2v_3 + v_1\otimes v_1v_2v_3 + v_3\otimes v_1^2v_2)\\
		&\qquad + c_3(v_2\otimes v_1^2v_3 + 2v_1\otimes v_2v_1v_3 + v_3\otimes v_1^2v_2)|.
	\end{align*}
	Therefore, $|a|\in \Ker(|\bar\Delta|)$ implies $2c_1 + c_2=0$, $c_2 + 2c_3 = 0$ and $c_1 + c_2 +  c_3= 0$. We have
	\[
		|a| = c_1|v_1^2v_2v_3  - 2v_1v_2v_1v_3 + v_2v_1^2v_3| = c_1 |[v_1,[v_1,v_2]]v_3|.
	\]
	Thus we have $|a|\in |HL(H)^{(3)}|$.
	\item  In multi-degree $(1,1,1,1,0,\dotsc,0)$, it is enough to consider the case $a = c_1v_1v_2v_3v_4  + c_2v_1v_3v_2v_4 + c_3v_2v_1v_3v_4 + c_4v_2v_3v_1v_4 + c_5v_3v_1v_2v_4 + c_6v_3v_2v_1v_4$. Then, 
	\begin{align*}
		0 &= |c_1(v_1\otimes v_2v_3v_4 + v_2\otimes v_1v_3v_4 + v_3\otimes v_1v_2v_4+ v_4\otimes v_1v_2v_3)\\
		&\qquad + c_2(v_1\otimes v_3v_2v_4 + v_2\otimes v_1v_3v_4 + v_3\otimes v_1v_2v_4+ v_4\otimes v_1v_3v_2)\\
		&\qquad + c_3(v_1\otimes v_2v_3v_4 + v_2\otimes v_1v_3v_4 + v_3\otimes v_2v_1v_4 + v_4\otimes v_2v_1v_3)\\
		&\qquad + c_4(v_1\otimes v_2v_3v_4 + v_2\otimes v_3v_1v_4 + v_3\otimes v_2v_1v_4 + v_4\otimes v_2v_3v_1)\\
		&\qquad + c_5(v_1\otimes v_3v_2v_4 + v_2\otimes v_3v_1v_4 + v_3\otimes v_1v_2v_4 + v_4\otimes v_3v_1v_2)\\
		&\qquad + c_6(v_1\otimes v_3v_2v_4 + v_2\otimes v_3v_1v_4 + v_3\otimes v_2v_1v_4 + v_4\otimes v_3v_2v_1)|.
	\end{align*}
	Therefore we have 
	\begin{align*}
		&c_1 + c_3 + c_4 = 0,\quad c_2 +c_5 + c_6 = 0,\\
		&c_1 + c_2 + c_3 =0 ,\quad c_4+c_5+c_6 = 0,\\
		&c_1+c_2+c_5=0,\quad c_3+c_4+c_6 = 0, \\
		&c_1+c_4+c_5=0,\textrm{ and }c_2+c_3+c_6=0.
	\end{align*}
	This boils down to $c_1=c_6$, $c_2=c_4$, $c_3=c_5$ and $c_1+c_2+c_3 = 0$. Finally, we have
	\begin{align*}
		|a| &= |c_1v_1v_2v_3v_4  + c_2v_1v_3v_2v_4 - (c_1+c_2)v_2v_1v_3v_4 + c_2v_2v_3v_1v_4 - (c_1+c_2)v_3v_1v_2v_4 + c_1v_3v_2v_1v_4|\\
		&= |c_1[[v_1,v_2],v_3]v_4 + c_2[[v_1,v_3],v_2]v_4|
	\end{align*}
	and we obtain $|a|\in |HL(H)^{(3)}|$.
\end{enumerate}
This completes the proof.\qed

\begin{remark}
Lemma \ref{lem:ckth} holds for an arbitrary $\mathbb{K}$-vector space $H$; the proof is identical.
\end{remark}

\subsection{Proof of Proposition \ref{prop:ckom}.} 

The kernel of the projection $|T(H)|^{\otimes 2} \to |T(H)_\omega|^{\otimes 2}$ is equal to $|\omid|\otimes |T(H)| \oplus |T(H)|\otimes |\omid|$. Since $|\omid|^{(\leq 2)} = 0$, the intersection with $\Imag|\bar\Delta| \subset |T(H)|\otimes H + H\otimes |T(H)|$ is zero in degree $d\leq 3$. Hence, the claim follows.

Now, we consider the degree 4 case. We modify the above multi-grading as follows. Take a symplectic basis $(x_i,y_i)_{1\leq i\leq g}$ of $H$ as before, and set $v_{2i-1} = x_i$ and $v_{2i} = y_i$ for $1\leq i\leq g$. The algebra $T(H)_\omega = T(H)/\omid$ is graded by the abelian group $D= \mathbb{Z}^{2g}/\langle e_1 + e_2 = e_3 + e_4 = \cdots = e_{2g-1} + e_{2g} \rangle$ where $(e_i)_{1\leq i\leq 2g}$ is the standard basis of $\mathbb{Z}^{2g}$ and $\deg(v_i)=e_i$. Indeed, $\omega$ is homogeneous with respect to the $D$-grading. We have a decomposition $T(H)_\omega = \bigoplus_{\mathbf{d}\in D} (T(H)_\omega)_\mathbf{d}$. Denoting the natural map by $p\colon (\mathbb{Z}_{\geq 0})^{2g} \to D$, the homogeneous part $(T(H)_\omega)_\mathbf{d}$ is non-zero if and only if $\mathbf{d}\in D^+ := p((\mathbb{Z}_{\geq 0})^{2g})$. The total degree map $|\cdot|\colon \mathbb{Z}^{2g}\to \mathbb{Z}$  descends to $|\cdot|\colon D \to \mathbb{Z}$.

Define the \textit{redundancy} $R\colon D^+ \to \mathbb{Z}$ by $R(\mathbf d) = \sum_{1\leq i\leq g} \mathrm{min}(\lambda_{2i-1},\lambda_{2i})$ with $\mathbf d = p\Big(\sum_i \lambda_i e_i\Big)$. This is well-defined, and we have $0\leq R(\mathbf d)\leq |\mathbf d|/2$ for any $\mathbf d$. Now that we set $|\mathbf d| = 4$, we only have to inspect the following cases:
\begin{enumerate}[(1)]\setcounter{enumi}{-1}
	\item If $R(\mathbf{d}) = 0$, we have $\#p^{-1}(\mathbf d) = 1$. Therefore, the natural map $T(H)_{p^{-1}(\mathbf d)} \twoheadrightarrow (T(H)_\omega)_\mathbf{d}$ is an isomorphism. The computation of $\Ker(|\bar\Delta_\omega|)$ is reduced to the case of $T(H)$.
	\item If $R(\mathbf{d}) = 1$, we have $\#p^{-1}(\mathbf d) = g$. The inverse image $p^{-1}(\mathbf d)$ is, up to permutations of pairs $(x_i,y_i)\leftrightarrow (x_j,y_j)$ and $(x_i,y_i)\leftrightarrow (y_i,-x_i)$, one of the followings:
	\begin{enumerate}[(i)]
		\item $\{e_{2i-1} + e_{2i} + 2e_{2g-1}:1\leq i\leq g\}$; or
		\item $\{e_{2i-1} + e_{2i} + e_{2g-3} + e_{2g-1}:1\leq i\leq g\}$.
	\end{enumerate}
	\item If $R(\mathbf{d}) = 2$, we have $\#p^{-1}(\mathbf d) = \frac{g(g+1)}{2}$. There is only one such $\mathbf d\in D^+$, and the inverse image is given by $p^{-1}(\mathbf d) = \{e_{2i-1} + e_{2i} + e_{2j-1} + e_{2j}:1\leq i\leq j\leq g\}$.
\end{enumerate}
Since the case $(0)$ is already done, we will deal with the cases $(1)$ and $(2)$ from now on. The case $g=1$ is easy, so we assume $g\geq 3$.

\noindent (1-i) It suffices to consider $a = \sum_{1\leq i<g} (c^i_1x_iy_ix_g^2 + c^i_2x_ix_gy_ix_g + c^i_3y_ix_ix_g^2) + c_1^g x_g^3y_g$. We may take $c^1_1 = - c_2^1/2$ by subtracting some scalar multiple of $\omega x_g^2$ from $a$. Then,
\begin{align*}
	0 &= \sum_{1\leq i<g} | c_1^i(x_i\otimes y_ix_g^2 + y_i\otimes x_ix_g^2 + 2x_g\otimes x_iy_ix_g)\\[-11pt]
	&\qquad\qquad  + c_2^i(x_i\otimes x_gy_ix_g + y_i\otimes x_ix_g^2 + x_g\otimes x_iy_ix_g + x_g\otimes x_ix_gy_i)\\
	&\qquad\qquad + c_3^i(x_i\otimes y_ix_g^2 + y_i\otimes x_ix_g^2 + 2x_g\otimes y_ix_ix_g)|\\
	&\qquad + c_1^g|3x_g\otimes x_g^2y_g + y_g\otimes x_g^3|. 
\end{align*}
We have to use two rewriting rules, but the only applicable term is $|y_{g-1}x_{g-1}x_g|$ in the seventh and tenth terms. By applying $\rho_2$ with $s=2$ and $t=1$, we have
\[
	|x_gy_{g-1}x_{g-1}| \mapsto |x_g(x_{g-1}y_{g-1} + \omega'')|.
\]
At this point, we can simply compare the coefficients by Theorem \ref{thm:basis}. Therefore, we have
\begin{align*}
	c_1^i + c_2^i + c_3^i = 0 &\quad \textrm{ for } 1\leq i< g,\\
	(c_2^i + 2c_3^i) - (c_2^{g-1} + 2c_3^{g-1}) = 0&\quad \textrm{ for } 1\leq i< g-1,\\
	(2c_1^i + c_2^i) + (c_2^{g-1} + 2c_3^{g-1}) = 0&\quad \textrm{ for } 1\leq i< g-1, \textrm{ and }\\
	c_1^g = 0&.
\end{align*}
Since we took $c^1_1 = - c_2^1/2$, we have $2c_1^i + c_2^i = c_2^i + 2c_3^i = 0$ for all $1\leq i<g$. The rest is the same as the case of $T(H)$.

\noindent (1-ii) It suffices to consider
\begin{align*}
	a &= \sum_{1\leq i<g-1} (c^i_1x_iy_ix_{g-1}x_g + c_2^ix_ix_{g-1}y_ix_g + c_3^iy_ix_ix_{g-1}x_g + c_4^iy_ix_{g-1}x_i x_g+ c_5^ix_{g-1}x_iy_ix_g + c_6^ix_{g-1}y_ix_ix_g)\\
	&\quad + (c_1^{g-1}x_{g-1}^2y_{g-1}x_g + c_2^{g-1} x_{g-1}y_{g-1}x_{g-1}x_g + c_3^{g-1}y_{g-1}x_{g-1}^2x_g)\\
	&\quad + (c_1^gx_{g-1}x_g^2y_g + c_2^g x_{g-1}x_gy_gx_g + c_3^g x_{g-1}y_gx_g^2).
\end{align*}
We have
\begin{align*}
	0 &= \sum_{1\leq i<g-1} |c_1^i(x_i\otimes y_ix_{g-1}x_g + y_i\otimes x_ix_{g-1}x_g + x_{g-1}\otimes x_iy_ix_g + x_g\otimes x_iy_ix_{g-1})\\[-11pt]
	&\qquad\qquad\quad + c_2^i(x_i\otimes x_{g-1}y_ix_g + y_i\otimes x_ix_{g-1}x_g + x_{g-1}\otimes x_iy_ix_g + x_g\otimes x_ix_{g-1}y_i)\\
	&\qquad\qquad\quad + c_3^i(x_i\otimes y_ix_{g-1}x_g + y_i\otimes x_ix_{g-1}x_g + x_{g-1}\otimes y_ix_ix_g + x_g\otimes y_ix_ix_{g-1})\\
	&\qquad\qquad\quad + c_4^i(x_i\otimes y_ix_{g-1}x_g + y_i\otimes x_{g-1}x_ix_g + x_{g-1}\otimes y_ix_ix_g + x_g\otimes y_ix_{g-1}x_i)\\
	&\qquad\qquad\quad + c_5^i(x_i\otimes x_{g-1}y_ix_g + y_i\otimes x_{g-1}x_ix_g + x_{g-1}\otimes x_iy_ix_g + x_g\otimes x_{g-1}x_iy_i)\\
	&\qquad\qquad\quad + c_6^i(x_i\otimes x_{g-1}y_ix_g + y_i\otimes x_{g-1}x_ix_g+ x_{g-1}\otimes y_ix_ix_g + x_g\otimes x_{g-1}y_ix_i)|\\
	&\qquad + c_1^{g-1}|2x_{g-1}\otimes x_{g-1}y_{g-1}x_g + y_{g-1}\otimes x_{g-1}^2x_g + x_g\otimes x_{g-1}^2y_{g-1}|\\
	&\qquad + c_2^{g-1}|x_{g-1}\otimes y_{g-1}x_{g-1}x_g + y_{g-1}\otimes x_{g-1}^2x_g + x_{g-1}\otimes x_{g-1}y_{g-1}x_g + x_g\otimes x_{g-1}^2y_{g-1}|\\
	&\qquad + c_3^{g-1}|2x_{g-1}\otimes y_{g-1}x_{g-1}x_g + y_{g-1}\otimes x_{g-1}^2x_g + x_g\otimes y_{g-1}x_{g-1}^2|\\
	&\qquad + c_1^g|x_{g-1}\otimes x_g^2y_g + 2x_g\otimes x_{g-1}x_gy_g + y_g\otimes x_{g-1}x_g^2|\\
	&\qquad + c_2^g|x_{g-1}\otimes x_g^2y_g + x_g\otimes x_{g-1}y_gx_g + y_g\otimes x_{g-1}x_g^2 + x_g\otimes x_{g-1}x_gy_g|\\
	&\qquad + c_3^g|x_{g-1}\otimes y_gx_g^2 + 2x_g\otimes x_{g-1}y_gx_g + y_g\otimes x_{g-1}x_g^2|.
\end{align*}
To apply the rewriting rules, we first search for the sequence $y_gx_g$, which appears in $c_2^g x_g\otimes x_{g-1}y_gx_g$ and $c_3^g\cdot 2x_g\otimes x_{g-1}y_gx_g$. The monomial $x_{g-1}y_gx_g$ is rewritten into $x_{g-1}(x_gy_g + \omega')$. Next, we search for the cyclic monomial $|x_gy_{g-1}x_{g-1}|$ or $|y_gy_{g-1}x_{g-1}|$. The latter does not appear since it has multi-degree $(0,\dotsc,0,1,1,0,1)$, which is impossible. The former appears in $c_2^{g-1}x_{g-1}\otimes y_{g-1}x_{g-1}x_g$ and $c_3^{g-1}\cdot 2x_{g-1}\otimes y_{g-1}x_{g-1}x_g$. After the rewriting, we can compare the coefficients; we have
\begin{align*}
	x_i : c_1^i + c_3^i + c_4^i = 0 &\quad \textrm{ for } 1\leq i< g-1,\\
	x_i : c_2^i + c_5^i + c_6^i = 0 &\quad \textrm{ for } 1\leq i< g-1,\\
	y_i : c_1^i + c_2^i + c_3^i = 0 &\quad \textrm{ for } 1\leq i< g-1,\\
	y_i : c_4^i + c_5^i + c_6^i = 0 &\quad \textrm{ for } 1\leq i< g-1,\\
	x_{g-1}: (c_1^i + c_2^i + c_5^i)  + (c_2^{g-1} + 2c_3^{g-1}) = 0 &\quad \textrm{ for } 1\leq i< g-1,\\
	x_{g-1}: (c_3^i + c_4^i + c_6^i)  - (c_2^{g-1} + 2c_3^{g-1}) = 0 &\quad \textrm{ for } 1\leq i< g-1,\\
	x_{g-1}: (c_1^{g-1} + c_2^{g-1} + c_3^{g-1})= 0 &,\\
	x_{g-1}: (c_1^{g} + c_2^{g} + c_3^{g})= 0 &,\\
	x_g: (c_1^i + c_4^i + c_5^i) + (c_2^g + 2c_3^g) = 0& \quad \textrm{ for } 1\leq i< g-1,\textrm{ and}\\
	x_g: (c_2^i + c_3^i + c_6^i) - (c_2^g + 2c_3^g) = 0& \quad \textrm{ for } 1\leq i< g-1.
\end{align*}
From this, we have $c_2^i = c_4^i$ for all $1\leq i< g-1$. Since we can freely take $c_1^1$ and $c_1^6$ by subtracting $\omega x_{g-1}x_g$ and $x_{g-1}\omega x_g$, respectively, we set $c_1^1 + c_4^1 + c_5^1 = c_2^1+c_3^1 + c_6^1 = 0$. Then, we have $c_2^g + 2c_3^g = 2c_1^c + c_2^g = 0$. Also we have $c_1^1 + c_2^1 + c_5^1 = 0$ since $c_2^1 = c_4^1$, so we get $c_2^{g-1}+2c_3^{g-1} = 0$ and $c_1^i + c_2^i + c_5^i = c_3^i + c_4^i + c_6^i = 0$ for all $i$. The rest is the same as the case of $T(H)$.

\noindent (2) It suffices to consider
\begin{align*}
	a &= \sum_{1\leq i \leq g} (c_1^ix_i^2y_i^2 + c_2^i x_iy_ix_iy_i)\\
	&\qquad + \sum_{1\leq i<j\leq g} (c_1^{ij}x_iy_ix_jy_j + c_2^{ij}x_iy_iy_jx_j + c_3^{ij}x_ix_jy_iy_j + c_4^{ij}x_ix_jy_jy_i + c_5^{ij}x_iy_jy_ix_j + c_6^{ij}x_iy_jx_jy_i).
\end{align*}
We have
\begin{align*}
	0 &= \sum_{1\leq i \leq g} |c_1^i(2x_i\otimes x_iy_i^2 + 2y_i\otimes x_i^2y_i) + c_2^i(2x_i\otimes x_iy_i^2 + 2y_i\otimes x_i^2y_i)|\\
	&+ \sum_{1\leq i<j\leq g} |c_1^{ij}(x_i\otimes y_ix_jy_j + y_i\otimes x_ix_jy_j + x_j\otimes x_iy_iy_j + y_j\otimes x_iy_ix_j)\\[-11pt]
	&\qquad\qquad\quad + c_2^{ij}(x_i\otimes y_iy_jx_j + y_i\otimes x_iy_jx_j + y_j\otimes x_iy_ix_j + x_j\otimes x_iy_iy_j)\\
	&\qquad\qquad\quad + c_3^{ij}(x_i\otimes x_jy_iy_j + x_j\otimes x_iy_iy_j + y_i\otimes x_ix_jy_j + y_j\otimes x_ix_jy_i)\\
	&\qquad\qquad\quad + c_4^{ij}(x_i\otimes x_jy_jy_i + x_j\otimes x_iy_jy_i + y_j\otimes x_ix_jy_i + y_i\otimes x_ix_jy_j)\\
	&\qquad\qquad\quad + c_5^{ij}(x_i\otimes y_jy_ix_j + y_j\otimes x_iy_ix_j + y_i\otimes x_iy_jx_j + x_j\otimes x_iy_jy_i)\\
	&\qquad\qquad\quad + c_6^{ij}(x_i\otimes y_jx_jy_i + y_j\otimes x_ix_jy_i + x_j\otimes x_iy_jy_i + y_i\otimes x_iy_jx_j)|\\
	&\xmapsto{\rho,\rho_2} \sum_{1\leq i \leq g} |c_1^i(2x_i\otimes x_iy_i^2 + 2y_i\otimes x_i^2y_i) + c_2^i(2x_i\otimes x_iy_i^2 + 2y_i\otimes x_i^2y_i)|\\
	&+ \sum_{1\leq i<j< g} |c_1^{ij}(x_i\otimes y_ix_jy_j + y_i\otimes x_ix_jy_j + x_j\otimes x_iy_iy_j + y_j\otimes x_iy_ix_j)\\[-11pt]
	&\qquad\qquad\quad + c_2^{ij}(x_i\otimes y_iy_jx_j + y_i\otimes x_iy_jx_j + y_j\otimes x_iy_ix_j + x_j\otimes x_iy_iy_j)\\
	&\qquad\qquad\quad + c_3^{ij}(x_i\otimes x_jy_iy_j + x_j\otimes x_iy_iy_j + y_i\otimes x_ix_jy_j + y_j\otimes x_ix_jy_i)\\
	&\qquad\qquad\quad + c_4^{ij}(x_i\otimes x_jy_jy_i + x_j\otimes x_iy_jy_i + y_j\otimes x_ix_jy_i + y_i\otimes x_ix_jy_j)\\
	&\qquad\qquad\quad + c_5^{ij}(x_i\otimes y_jy_ix_j + y_j\otimes x_iy_ix_j + y_i\otimes x_iy_jx_j + x_j\otimes x_iy_jy_i)\\
	&\qquad\qquad\quad + c_6^{ij}(x_i\otimes y_jx_jy_i + y_j\otimes x_ix_jy_i + x_j\otimes x_iy_jy_i + y_i\otimes x_iy_jx_j)|\\
	&+ \sum_{1\leq i<g-1} |c_1^{ig}(x_i\otimes y_ix_gy_g + y_i\otimes x_ix_gy_g + x_g\otimes x_iy_iy_g + y_g\otimes x_iy_ix_g)\\[-11pt]
	&\qquad\qquad\quad + c_2^{ig}(x_i\otimes y_i(x_gy_g + \omega') + y_i\otimes x_i(x_gy_g + \omega') + y_g\otimes x_iy_ix_g + x_g\otimes x_iy_iy_g)\\
	&\qquad\qquad\quad + c_3^{ig}(x_i\otimes y_i(x_gy_g + \omega') + x_g\otimes x_iy_iy_g + y_i\otimes x_ix_gy_g + y_g\otimes x_ix_gy_i)\\
	&\qquad\qquad\quad + c_4^{ig}(x_i\otimes x_gy_gy_i + x_g\otimes x_iy_gy_i + y_g\otimes x_ix_gy_i + y_i\otimes x_ix_gy_g)\\
	&\qquad\qquad\quad + c_5^{ig}(x_i\otimes y_gy_ix_g + y_g\otimes x_iy_ix_g + y_i\otimes x_i(x_gy_g + \omega') + x_g\otimes x_iy_gy_i)\\
	&\qquad\qquad\quad + c_6^{ig}(x_i\otimes (x_gy_g + \omega')y_i + y_g\otimes x_ix_gy_i + x_g\otimes x_iy_gy_i + y_i\otimes x_i(x_gy_g + \omega'))|\\
	&+ |c_1^{g-1,g}(x_{g-1}\otimes y_{g-1}x_gy_g + y_{g-1}\otimes x_{g-1}x_gy_g + x_g\otimes x_{g-1}y_{g-1}y_g + y_g\otimes x_{g-1}y_{g-1}x_g)\\
	&\quad + c_2^{g-1,g}(x_{g-1}\otimes y_{g-1}(x_gy_g + \omega') + y_{g-1}\otimes x_{g-1}(x_gy_g + \omega') + y_g\otimes x_{g-1}y_{g-1}x_g + x_g\otimes x_{g-1}y_{g-1}y_g)\\
	&\quad + c_3^{g-1,g}(x_{g-1}\otimes y_{g-1}(x_gy_g + \omega') + x_g\otimes x_{g-1}y_{g-1}y_g + y_{g-1}\otimes x_{g-1}x_gy_g + y_g\otimes x_g(x_{g-1}y_{g-1} + \omega''))\\
	&\quad + c_4^{g-1,g}(x_{g-1}\otimes x_gy_gy_{g-1} + x_g\otimes y_g(x_{g-1}y_{g-1} + \omega'') + y_g\otimes x_g(x_{g-1}y_{g-1} + \omega'') + y_{g-1}\otimes x_{g-1}x_gy_g)\\
	&\quad + c_5^{g-1,g}(x_{g-1}\otimes y_gy_{g-1}x_g + y_g\otimes x_{g-1}y_{g-1}x_g + y_{g-1}\otimes x_{g-1}(x_gy_g + \omega') + x_g\otimes y_g(x_{g-1}y_{g-1} + \omega''))\\
	&\quad + c_6^{g-1,g}(x_{g-1}\otimes (x_gy_g + \omega')y_{g-1} + y_g\otimes x_g(x_{g-1}y_{g-1} + \omega'')\\
	&\qquad\qquad\quad + x_g\otimes y_g(x_{g-1}y_{g-1} + \omega'') + y_{g-1}\otimes x_{g-1}(x_gy_g + \omega'))|.
\end{align*}
Now, we compare the coefficients. We obtain
\begin{align}
	x_i\otimes x_iy_i^2&: 2c_1^i + 2c_2^i = 0\,(1\leq i\leq g),\\
	x_i\otimes y_ix_jy_j&: (c_1^{ij} + c_4^{ij} + c_5^{ij}) + (c_2^{ig} + c_3^{ig} + c_6^{ig}) = 0 \,(1\leq i<j<g),\label{eq:ijx1}\\
	x_i\otimes y_iy_jx_j&: (c_2^{ij} + c_3^{ij} + c_6^{ij}) - (c_2^{ig} + c_3^{ig} + c_6^{ig}) = 0 \,(1\leq i<j<g),\label{eq:ijx2}\\
	y_i\otimes x_ix_jy_j&: (c_1^{ij} + c_3^{ij} + c_4^{ij}) + (c_2^{ig} + c_5^{ig} + c_6^{ig}) = 0 \,(1\leq i<j<g),\label{eq:ijy1}\\
	x_j\otimes x_iy_iy_j&: (c_1^{ij} + c_2^{ij} + c_3^{ij}) + (c_2^{jg} + c_3^{jg} + c_6^{jg}) = 0 \,(1\leq i<j<g),\label{eq:jix1}\\
	x_j\otimes x_iy_jy_i&: (c_4^{ij} + c_5^{ij} + c_6^{ij}) - (c_2^{jg} + c_3^{jg} + c_6^{jg}) = 0 \,(1\leq i<j<g),\label{eq:jix2}\\
	y_j\otimes x_iy_ix_j&: (c_1^{ij} + c_2^{ij} + c_5^{ij}) + (c_2^{jg} + c_5^{jg} + c_6^{jg}) = 0 \,(1\leq i<j<g),\label{eq:jiy1}\\
	x_i\otimes y_ix_gy_g&: c_1^{ig} + c_2^{ig} + c_3^{ig} + c_4^{ig} + c_5^{ig} + c_6^{ig} = 0\,(1\leq i<g),\label{eq:igx0}\\
	x_g\otimes x_iy_iy_g&: (c_1^{ig} + c_2^{ig} + c_3^{ig}) + (c_4^{g-1,g} + c_5^{g-1,g} + c_6^{g-1,g}) = 0\,(1\leq i<g-1),\label{eq:igx1}\\
	x_g\otimes y_ix_iy_g&: (c_4^{ig} + c_5^{ig} + c_6^{ig}) - (c_4^{g-1,g} + c_5^{g-1,g} + c_6^{g-1,g}) = 0\,(1\leq i<g-1), \textrm{ and}\label{eq:igx2}\\
	y_g\otimes x_iy_ix_g&: (c_1^{ig} + c_2^{ig} + c_5^{ig}) + (c_3^{g-1,g} + c_4^{g-1,g} + c_6^{g-1,g}) = 0\,(1\leq i<g-1).\label{eq:igy1}
\end{align}
Equations \eqref{eq:igx1}-\eqref{eq:igy1} holds also for $i=g-1$: \eqref{eq:igx1} and \eqref{eq:igy1} follows from \eqref{eq:igx0}, and \eqref{eq:igx2} is a tautology. Now put $b^{ij} = c_3^{ij} - c_5^{ij}$. From \eqref{eq:ijx1} and \eqref{eq:ijy1}, we have $b^{ij} - b^{ig} = 0$ for $1\leq i<j<g$. From \eqref{eq:jix1} and \eqref{eq:jiy1}, we have $b^{ij} + b^{jg} = 0$ for $1\leq i<j<g$. From \eqref{eq:igx1} and \eqref{eq:igy1}, we have $b^{ig} - b^{g-1,g} = 0$ for $1\leq i < g$. Now that we have assumed $g\geq 3$, we have
\[
	-b^{2g} = b^{12} = b^{1g} = b^{g-1,g} = b^{2g},
\]
which implies $b^{2g}=0$ since $\mathrm{char}(\mathbb{K})\neq 2$. Hence we have $b^{ij} = 0$ for $1\leq i<j\leq g$, which says $c_3^{ij} = c_5^{ij}$ for $1\leq i<j\leq g$. By subtracting a scalar multiple of $x_iy_i\omega$ and $y_ix_i\omega$ for $1\leq i < g$, we can freely set $c_1^{ig}$ and $c_6^{ig}$ to any value. We set $c_1^{ig} = c_6^{ig} = -c_2^{ig} - c_3^{ig}$ for $1\leq i<g$. Then, from \eqref{eq:ijx1}-\eqref{eq:jiy1}, we have
\[
	c_1^{ij} + c_4^{ij} + c_5^{ij} = c_2^{ij} + c_3^{ij} + c_6^{ij} =  c_1^{ij} + c_3^{ij} + c_4^{ij} = c_1^{ij} + c_2^{ij} + c_3^{ij} = c_4^{ij} + c_5^{ij} + c_6^{ij} = c_1^{ij} + c_2^{ij} + c_5^{ij} = 0.
\]
for $1\leq i<j<g$. This gives $c_1^{ij} = c_6^{ij}$ and $c_2^{ij} = c_4^{ij}$ for $1\leq i<j<g$. Furthermore, since we set $c_1^{ig} + c_2^{ig} + c_3^{ig} = 0$, combining with \eqref{eq:igx0}, we have $c_4^{ig} + c_5^{ig} + c_6^{ig} = 0$ for any $i$. Then, by \eqref{eq:igx1} and \eqref{eq:igx2}, we have
\[
	c_1^{ig} + c_2^{ig} + c_3^{ig} = c_4^{ig} + c_5^{ig} + c_6^{ig} = 0
\]
and, together with $c_1^{ig} = c_6^{ig}$ and $c_3^{ig} = c_5^{ig}$ we already know, we have $c_2^{ig} = c_4^{ig}$ for all $i$. In conclusion, we have $c_1^i + c_2^i = 0$ for all $i$ and $c_1^{ij} = c_6^{ij}$, $c_2^{ij} = c_4^{ij}$, $c_3^{ij} = c_5^{ij}$, $c_1^{ij} + c_2^{ij} + c_3^{ij} = 0$ for $1\leq i<j\leq g$. The rest is the same as before.\qed

\begin{remark}
The author conjectures that $\Ker(|\bar\Delta|)^{(5)} = |HL^{(4)}| \oplus (\wedge^5H)$, whose rigorous proof could be given by a straightforward generalisation of the above procedure. The equality is already verified by a computer-assisted calculation.
\end{remark}

We conclude this section with a conjecture.

\begin{conjecture}\label{conj:surj}
The map $\Ker(|\bar\Delta|)^{(d)} \to \Ker(|\bar\Delta_\omega|)^{(d)}$ induced from the natural surjection $|T(H)| \twoheadrightarrow |T(H)_\omega|$ is also a surjection except for $(g,d)= (2,4)$.\\
\end{conjecture}

\small
\bibliographystyle{alphaurl}
\bibliography{kvcs.bib}

\end{document}